**Modeling and computation of an integral operator Riccati equation for an infinite-dimensional stochastic differential equation governing streamflow discharge**

Hidekazu Yoshioka[a, *], Motoh Tsujimura[b], Tomohiro Tanaka[c], Yumi Yoshioka[a], Ayumi Hashiguchi[a]

[a]Faculty of Life and Environmental Sciences, Shimane University, Nishikawatsu-cho 1060, Matsue 690-8504, Japan. E-mail: yoshih@life.shimane-u.ac.jp (H. Yoshioka), yyoshioka@life.shimane-u.ac.jp (Y. Yoshioka), hashiguchi@life.shimane-u.ac.jp (A. Hashiguchi)

[b]Graduate School of Commerce, Doshisha University, Karasuma-Higashi-iru, Imadegawa-dori, Kamigyo-ku, Kyoto 602-8580, Japan. E-mail: mtsujimu@mail.doshisha.ac.jp

[c]Graduate School of Engineering, Kyoto University, Kyotodaigaku-katsura, Nishikyo-ku, Kyoto, Kyoto 615-8540, Japan. E-mail: tanaka.tomohiro.7c@kyoto-u.ac.jp

* Corresponding author: Assistant Professor, Shimane University, Nishikawatsu-cho 1060, Matsue 690-8504, Japan. Tel: +81-852-32-6541, E-mail: yoshih@life.shimane-u.ac.jp

**Contributions**

Hidekazu Yoshioka: Authorization, Data acquisition, Data curation, Field survey, Formal analysis, Numerical computation, Writing and editing, Funding acquisition

Motoh Tsujimura: Formal analysis, Writing and editing

Tomohiro Tanaka: Numerical computation, Writing and editing

Yumi Yoshioka: Data acquisition, Field survey, Writing and editing

Ayumi Hashiguchi: Data acquisition, Field survey, Writing and editing

**Abstract** We propose a linear-quadratic (LQ) control problem of streamflow discharge by optimizing an infinite-dimensional jump-driven stochastic differential equation (SDE). Our SDE is a superposition of Ornstein–Uhlenbeck processes (supOU process), generating a sub-exponential autocorrelation function observed in actual data. The integral operator Riccati equation is heuristically derived to determine the optimal control of the infinite-dimensional system. In addition, its finite-dimensional version is derived with a discretized distribution of the reversion speed and computed by a finite difference scheme. The optimality of the Riccati equation is analyzed by a verification argument. The supOU process is parameterized based on the actual data of a perennial river. The convergence of the numerical scheme is analyzed through computational experiments. Finally, we demonstrate the application of the proposed model to realistic problems along with the Kolmogorov backward equation for the performance evaluation of controls.

**Classification Codes** 93E20, 49L12, 60H10, 45K05, 65Y04

**Abbreviations**

| | |
|---|---|
| ACF | Autocorrelation function |
| HJB | Hamilton–Jacobi–Bellman |
| KBE | Kolmogorov backward equation |
| LQ | Linear-quadratic |
| MC | Monte-Carlo |
| OU | Ornstein–Uhlenbeck |
| PDF | Probability density function |
| SDE | Stochastic differential equation |
| supOU | superposition of Ornstein–Uhlenbeck |

# 1 Introduction

## 1.1 Research background

Rivers are essential elements in the transport of water and nutrients (Revel et al., 2021; Drake et al., 2021). In addition, they support aquatic ecosystems by supplying habitats for a variety of organisms (Kaiser et al., 2020; Watz et al., 2021). Therefore, designing streamflow regimes to support the river environment has been a hot research topic (Guan et al., 2021; Wu et al., 2022; Tonkin et al., 2021). Optimization approaches have long been employed for streamflow design. In particular, models based on stochastic processes where hydrological variables such as rainfall and discharge vary stochastically have been employed because of their simplicity and flexibility (Farzaneh et al., 2021; Giuliani et al., 2021; Yoshioka, 2022).

Streamflow discharge, hereafter called discharge, is a stochastic process that plays a crucial role in modeling river environments. A classical tool for stochastic streamflow modeling is the jump-driven stochastic differential equation (SDE), which governs the occurrence and recession of floods (Botter et al., 2010). Ornstein–Uhlenbeck (OU) processes have been employed for hydrological analysis, including reproducing flow duration curves (Santos et al., 2018), biogeochemical dynamics (Musolff et al., 2017), hydrological similarity analysis (Betterle et al., 2019), stability analysis of plant uprooting (Calvani et al., 2019), and water quality analysis (Yoshioka and Yoshioka, 2021).

A drawback of OU processes is the exponentially decaying autocorrelation functions (ACFs) ($\exp(-c\tau)$, where $c > 0$, and $\tau$ is the time lag), which are unsatisfactory when the discharge time series exhibit a sub-exponential and hence slower power-type decay ($\tau^{-c}$ for large $\tau$) (Habib, 2020; Al Sawaf et al., 2017; Mohamed et al., 2017; Mihailović et al., 2019). The sub-exponential decay implies a more persistent memory than those derived by OU processes. Because the ACF characterizes temporal timescales, the sub-exponential decay must be considered in stochastic process models. From a practical viewpoint, the flexibility of OU processes should be kept as well.

We apply a superposition of OU processes, a supOU process (Barndorff-Nielsen and Stelzer, 2011), to discharge time series. The superposition is conducted for the reversion speed to allow for multiple timescales responsible for the sub-exponential decay. This approach has recently been applied to efficiently model a mountainous river environment (Yoshioka, 2022). Botter (2010) considered an analogous jump-driven OU process whose recession rate changes probabilistically at each jump event, focusing solely on marked Poisson jumps. The OU-type models driven by tempered stable processes that involve infinite activities have been found to better describe streamflow dynamics (Yoshioka and Tsujimura, 2022; Yoshioka and Yoshioka, 2021). supOU and related processes have also been applied to time series that have long memory in finance and related research areas (Barndorff-Nielsen and Stelzer, 2013; Moser and Stelzer, 20111); however, their application to environmental optimization problems has yet to be addressed.

In contrast to the conventional Markovian OU processes, supOU processes are non-Markovian processes that can reproduce sub-exponential ACFs (Stelzer et al., 2015). The non-Markovian nature comes from the probabilistically distributed reversion speed, which by contrast is a constant in classical OU processes. This non-Markovian property hinders us from applying the standard machinery of stochastic control (e.g., Pham (2009)) to supOU processes. Such an issue is ubiquitous in generic non-Markovian models including Volterra processes (Barndorff-Nielsen et al., 2018). Unfortunately, explicit solutions have been found for only limited cases (Hainaut, 2022; Han and Wong, 2021; Han et al., 2021).

Recently, progress has been made in the analysis and computation of infinite-dimensional optimality equations in stochastic control. Alfonsi and Schied (2013) mathematically and numerically analyzed the well-posedness and optimality of the infinite-dimensional Riccati equation of a singular control type. Abi Jaber et al. (2021a-b) analyzed the finite-dimensional approximation errors of the infinite-dimensional Riccati equations of the linear-quadratic (LQ) control problems of stochastic Volterra processes. Their methodologies were based on the Markovian lift where a non-Markovian process was reformulated as an infinite sum of Markovian processes to which the Markovian stochastic calculus applies (Abi Jaber, 2019). Similar techniques have been used in finance and economics (Bäuerle and Desmettre, 2020; Jusselin, 2021). Infinite-dimensional Riccati equations in an infinite horizon have been studied to control partial differential equations (Da Prato and Ichikawa, 1988; Bay et al., 2011), while their engineering applications have remained rare. Infinite-dimensional control problems have already appeared in environmental and ecological management under state uncertainty where an infinite-dimensional optimality equation is reduced to a finite-dimensional equation with specific uncertainty forms (Kling et al., 2017; Sloggy et al., 2020). These problems differ from LQ ones, but both share the philosophy of reducing an infinite-dimensional optimality equation to a finite-dimensional and computable equation by restricting the problem class.

### 1.2 Objective and contribution

We present a mathematical model for the temporal evolution of discharge based on a supOU process, the formulation and analysis of its LQ control problem, and its real-world application. Thus, we aimed to develop a non-Markovian streamflow discharge model. A long-run LQ control problem is formulated along with its approximation by a Markovian lift to minimize the deviation between the actual and seasonally varying, targeted discharges by adding/subtracting water. While LQ controls have been applied to stormwater management (Wong and Kerkez, 2018) and irrigation systems (Zhong et al., 2020; Conde et al., 2021), in this work, we introduce the use of supOU processes to address control problems in river management.

We explore periodic solutions to the corresponding infinite-dimensional integral operator Riccati equation both theoretically and numerically. Our control problem and the Riccati equation are therefore different from those used in stationary processes as they have time-dependent parameters (Bardi and Priuli, 2014; Priuli, 2015; Guatteri and Masiero, 2009; Bian and Jiang, 2019; Mei et al., 2021). For a reduced finite-dimensional case, the well-posedness of the problem follows from similarities between our and existing Riccati equations. By contrast, for the infinite-dimensional case, because of the high irregularity of noise, we heuristically derive the Riccati equation. We discuss that the Riccati equation is an optimality equation of the control problem under certain assumptions. The heuristics are supplemented by numerical experiments. The convergence of the finite-dimensional Riccati equation to the infinite-dimensional one is numerically analyzed.

Using a Lagrangian multiplier technique (e.g., Dai et al, 2010; Pfeiffer, 2018; Euch et al., 2021) and the Kolmogorov backward equation (KBE), we analyze the cost–performance relationship and stabilizability of controls. Indeed, the stabilizability of our control is not trivial because of the lack of classical dissipativity. The KBE provides a methodology to evaluate the stability of the controlled dynamics. Thus, we avoid directly handling infinite-dimensional SDEs. This KBE is used to analyze the performance of the optimal control.

To apply the proposed model to realistic cases, we identify the supOU processes in a perennial river. The supOU process at each observation station is identified from the available data. The identified model is finally applied to an LQ control problem of discharge to effectively suppress the growth of nuisance benthic algae (Huang et al., 2021; Loire et al., 2021; Yoshioka and Tsujimura, 2020). Our contribution offers a unique approach to solving this environmental problem through the optimal control of a non-Markovian stochastic process model.

### 1.3 Outline of this work

In **Section 2**, we firstly formulate supOU processes and then a LQ long-term control problem. Secondly, we present a finite-dimensional version of a supOU process, which plays a vital role in numerical computation. The consistency between the finite- and infinite-dimensional processes is analyzed through characteristic functions. In **Section 3**, the Hamilton–Jacobi–Bellman (HJB) equation of the finite-dimensional case is obtained by dynamic programming, and the associated Riccati equation is obtained by a quadratic ansatz. The corresponding infinite-dimensional problem is formulated as well. We also evaluate the total cost of each optimal control. For this purpose, the KBEs for the finite- and infinite-dimensional cases are presented to compute the cost and performance. The optimality of the control derived from the Riccati equation is theoretically verified as well. In **Section 4**, our model is applied to a perennial river along with numerical experiments to analyze the convergence of numerical solutions. In **Section 5**, we provide the conclusions and implications of the study.

**Appendix A** provides supporting data on emulating a discharge time series at an observation station of our study site. **Appendix B** presents jump processes alternative to what we use in this study. **Appendix C** contains proofs of **Propositions 1-3**. **Appendix D** discusses the control problem under a simplified setting.

## 2 Mathematical model

### 2.1 supOU process

A supOU process is a superposition of OU processes with respect to the reversion speed (e.g., Barndorff-Nielsen and Stelzer, 2011; Griffin, 2010; Barrera and Pardo, 2020). Our supOU process

$X = (X_t)_{t\geq 0}$ is a stationary process given by

$$X_t = \underline{X} + \int_{-\infty}^{t}\int_0^{+\infty}\int_0^{+\infty} e^{-\lambda(t-s)} z\mu(\mathrm{d}z, \mathrm{d}\lambda, \mathrm{d}s),\quad t \in \mathbb{R}, \tag{1}$$

where $\underline{X} \geq 0$ is a constant, $\mu$ is a (double-sided in time) random Poisson measure on $(0,+\infty)\times(0,+\infty)\times\mathbb{R}$ (a product space of (jump size) × (reversion speed) × (time)) with the compensator $\nu(\mathrm{d}z)\pi(\mathrm{d}\lambda)\mathrm{d}s$ containing a probability measure $\pi$ on $(0,+\infty)$ and the background driving Lévy measure $\nu$ (e.g., Peszat and Zabczyk, 2007; Barndorff-Nielsen and Stelzer, 2011). We assume that $\nu$ is a positive-jump subordinator type satisfying

$$\int_0^{+\infty} z^2\nu(\mathrm{d}z),\ \int_0^{+\infty}\min\{1,z\}\nu(\mathrm{d}z) < +\infty, \tag{2}$$

generating a Lévy process with bounded variation. For $\pi$, we assume

$$\int_0^{+\infty}\pi(\mathrm{d}\lambda) = 1,\ \int_0^{+\infty}\frac{1}{\lambda}\pi(\mathrm{d}\lambda) < +\infty, \tag{3}$$

meaning the usual normalization and boundedness of the reciprocal moment. Under (2)-(3), the supOU process (1) exists globally in time as an infinitely divisible stationary process (e.g., Theorem 3.1 of Barndorff-Nielsen and Stelzer 2011; Theorem 2.2 of Stelzer et al., 2015). The supOU process (1) is adapted to the filtration generated by $\mu$ (Theorem 3.12 of Barndorff-Nielsen and Stelzer, 2011). This filtration is denoted with the sigma algebra $\mathcal{F}_t$ at time $t$ as $\mathcal{F} = (\mathcal{F}_t)_{t\geq 0}$.

We use the following $\nu$ and $\pi$ to focus on the application in **Section 4**:

$$\nu(\mathrm{d}z) = a_\nu \frac{e^{-b_\nu z^{p_\nu}}}{z^{1+\alpha_\nu}}\mathrm{d}z,\quad z > 0 \tag{4}$$

and

$$\pi(\mathrm{d}\lambda) = \frac{1}{\Gamma(\alpha_\pi)}\frac{\lambda^{\alpha_\pi - 1}}{B_\pi^{\alpha_\pi}} e^{-\frac{\lambda}{B_\pi}}\mathrm{d}\lambda,\quad \lambda > 0, \tag{5}$$

where $\Gamma(\cdot)$ is a Gamma function. In (4), $a_\nu > 0,\ b_\nu > 0,\ \alpha_\nu < 1,\ p_\nu > 0$ are parameters with which $\nu$ becomes a generalized tempered stable type as a flexible subordinator (Grabchak, 2021). It reduces to a classical tempered stable one if $p_\nu = 1$. Jumps generated by $\nu$ have finite (resp., infinite) activities if $\alpha_\nu < 0$ (resp., $0 \leq \alpha_\nu < 1$); thus, (4) reproduces a wide variety of jump processes by tuning $\alpha_\nu$. In (5), $\alpha_\pi > 1,\ B_\pi > 0$ are shape and scaling parameters. This $\pi$ implies the existence of reversion speeds ranging from 0 to $+\infty$ and thus multiple timescales in each sample path of $X$. Phenomenologically, each jump associates a reversion speed, leading to the existence of a variety of reversion speeds in each sample path (Stelzer et al., 2015). This multiscale property is essential for sub-exponential ACFs.

The measures in (4) and (5) allow for explicitly obtaining moments and the ACF of $X$.

$$M_k \equiv \int_0^{+\infty} z^k\nu(\mathrm{d}z) = \frac{a_\nu}{p_\nu}(b_\nu)^{\frac{\alpha_\nu - k}{p_\nu}}\Gamma\left(\frac{k-\alpha_\nu}{p_\nu}\right)\ (k\in\mathbb{N}),\quad R \equiv \int_0^{+\infty}\frac{1}{\lambda}\pi(\mathrm{d}\lambda) = \frac{1}{B_\pi(\alpha_\pi - 1)}. \tag{6}$$

Based on the work by Yoshioka (2022), we describe the stationary moments as

$$\text{Ave} = \mathbb{E}[X_t] = \underline{X} + RM_1 , \tag{7}$$

$$\text{Var} = \mathbb{E}\left[\left(X_t - \mathbb{E}[X_t]\right)^2\right] = \frac{RM_2}{2} , \tag{8}$$

$$\text{Skew} = \frac{\mathbb{E}\left[\left(X_t - \mathbb{E}[X_t]\right)^3\right]}{\text{Var}^{3/2}} = \frac{RM_3}{3\text{Var}^{3/2}} , \tag{9}$$

$$\text{Kurt} = \frac{\mathbb{E}\left[\left(X_t - \mathbb{E}[X_t]\right)^4\right]}{\text{Var}^2} - 3 = \frac{RM_4}{4\text{Var}^2} , \tag{10}$$

and the ACF with a time lag $\tau \geq 0$

$$\text{ACF}(\tau) = \left(\frac{1}{1 + B_\pi \tau}\right)^{\alpha_\pi - 1} . \tag{11}$$

The statistical quantities (7)–(11) are used for identifying model parameters.

## 2.2 Truncated supOU process

### 2.2.1 Formulation

The truncated, finite-dimensional supOU process is introduced because it serves as the foundation of our numerical computation of a finite-dimensional Riccati equation. Set the truncated supOU process

$$X_{n,t} = \underline{X} + \sum_{i=1}^{n} \int_{-\infty}^{t} e^{-\lambda_i (t-s)} \mathrm{d}L_s^{(i)} , \quad t \in \mathbb{R} \tag{12}$$

with a non-negative sequence $\{c_i\}_{i=1,2,3,...}$ with $\sum_{i=1}^{\infty} c_i = 1$, and $X_t^{(i)}$ is the OU process at $t$ with the recession rate $\lambda_i$ and the pure-jump Lévy process $L^{(i)}$ having the measure $c_i v$ with each $L_s^{(i)}$ ($i = 1, 2, 3, ...$) being mutually independent (Fasen and Klüppelberg, 2007).

We show that the truncated supOU process (12) approximates the original one (1) with appropriate $\{c_i, \lambda_i\}_{i=1,2,3,...}$ in the sense of law. For each $n$, we choose a strictly increasing sequence $\eta_n = \{\eta_{n,i}\}_{1 \leq i \leq n}$ with $\eta_{n,0} = 0$. We set the uniform mesh to discretize the space of reversion speed

$$\eta_{n,i} = \bar{\eta} \frac{i}{n^\beta} , \quad 1 \leq i \leq n \tag{13}$$

with $\beta \in (0,1)$, $\bar{\eta} > 0$. Furthermore,

$$c_i = \int_{\eta_{n,i-1}}^{\eta_{n,i}} \pi(\mathrm{d}\lambda) , \quad \lambda_i = \frac{1}{c_i} \int_{\eta_{n,i-1}}^{\eta_{n,i}} \pi(\mathrm{d}\lambda) \in (\eta_{n,i-1}, \eta_{n,i}) , \quad 1 \leq i \leq n \tag{14}$$

and the discrete measure

$$\pi_n(\mathrm{d}\lambda) = \sum_{i=1}^{n} c_i \delta_{\{\lambda = \lambda_i\}} \quad \text{with} \quad c_{n+1} = 1 - \sum_{i=1}^{n} c_i \quad \text{(introduced for brevity) and} \quad \lambda_{n+1} = +\infty . \tag{15}$$

### 2.2.2 Consistency

For simplicity, set $\underline{X}=0$ in this sub-section without any loss of generality. Set the imaginary unit as $\mathrm{i}$ ($\mathrm{i}^2=-1$). The principal branch of the characteristic function $\varphi_{X_t}(u)$ ($u\in\mathbb{R}$) of the supOU process (1) at time $t\in\mathbb{R}$ is (Proposition 2.4 of Barndorff-Nielsen and Stelzer, 2011)

$$\ln\varphi_{X_t}(u)=\ln\mathbb{E}\left[\exp(\mathrm{i}uX_t)\right]=\int_0^{+\infty}(\exp(\mathrm{i}uz)-1)\bar{v}(\mathrm{d}z) \tag{16}$$

with

$$\bar{v}(Z)=\int_0^{+\infty}\int_{-\infty}^{t}\int_0^{+\infty}\chi_Z\left(e^{-\lambda(t-s)}z\right)v(\mathrm{d}z)\mathrm{d}s\pi(\mathrm{d}\lambda) \tag{17}$$

for any Borel measurable set $Z$ of $(0,\infty)$ with an indicator function $\chi_Z(\cdot)$ such that $\chi_Z(y)=1$ if $y\in Z$ and $\chi_Z(y)=0$ otherwise. By (17), (16) becomes

$$\begin{aligned}\ln\varphi_{X_t}(u)&=\int_0^{+\infty}(\exp(\mathrm{i}uz)-1)\bar{v}(\mathrm{d}z)\\&=\int_0^{+\infty}\int_{-\infty}^{t}\int_0^{+\infty}\left(\exp\left(\mathrm{i}ue^{-\lambda(t-s)}z\right)-1\right)\chi_{(0,+\infty)}\left(e^{-\lambda(t-s)}z\right)v(\mathrm{d}z)\mathrm{d}s\pi(\mathrm{d}\lambda).\\&=\int_0^{+\infty}\int_{-\infty}^{t}\int_0^{+\infty}\left(\exp\left(\mathrm{i}ue^{-\lambda(t-s)}z\right)-1\right)v(\mathrm{d}z)\mathrm{d}s\pi(\mathrm{d}\lambda)\end{aligned} \tag{18}$$

Similarly, the characteristic function $\varphi_{n,X_{n,t}}(\cdot)$ of the truncated supOU process at $t\in\mathbb{R}$ is

$$\ln\varphi_{n,X_{n,t}}(u)=\int_0^{+\infty}\int_{-\infty}^{t}\int_0^{+\infty}\left(\exp\left(\mathrm{i}ue^{-\lambda(t-s)}z\right)-1\right)v(\mathrm{d}z)\mathrm{d}s\pi_n(\mathrm{d}\lambda). \tag{19}$$

The main result here is the following **Proposition 1** stating a consistency of the finite-dimensional supOU process in the sense of characteristics functions. For brevity, set

$$\bar{G}(n)=\frac{1}{\eta_{n,n}}\int_{\eta_{n,n}}^{+\infty}\pi(\mathrm{d}\lambda)+\sum_{i=1}^{n}\frac{(\eta_{n,i}-\eta_{n,i-1})^2}{\lambda_i},\quad n\in\mathbb{N}. \tag{20}$$

The proof of **Proposition 1** is placed in **Appendix C.1**. As shown in this appendix, $\bar{G}(n)$ corresponds to the convergence rate of $\left|\ln\varphi_{X_t}(u)-\ln\varphi_{n,X_{n,t}}(u)\right|$.

***Proposition 1***

*Assume* $\alpha_\pi>2$ *and* $\bar{G}(n)\to 0$ *as* $n\to+\infty$. *Then,* $\left|\varphi_{X_t}(u)-\varphi_{n,X_{n,t}}(u)\right|\to 0$ *as* $n\to+\infty$ *pointwise in* $\mathbb{R}^2$. *Hence,* $X_{n,t}$ *converges to* $X_t$ *in the sense of law as* $n\to+\infty$.

The assumption $\alpha_\pi>2$ has a technical origin but is satisfied in our application in **Section 4**. We expect $\lambda_i\approx\frac{1}{2}(\eta_{n,i}+\eta_{n,i-1})$ for $n>>1$ because $\lambda^{\alpha_\pi-1}e^{-\frac{\lambda}{B_\pi}}$ ($\lambda\geq 0$) is continuously differentiable and Lipschitz continuous, suggesting $\lim_{n\to+\infty}\sum_{i=1}^{n}\frac{\left(\eta_{n,i}-\eta_{n,i-1}\right)^2}{\lambda_i}=\lim_{n\to+\infty}\frac{\overline{\eta}}{n^{\beta}}\sum_{i=1}^{n}\frac{1}{i-1/2}=\lim_{n\to+\infty}\frac{\overline{\eta}}{n^{\beta}}\ln n\to+0$ hence the assumption of **Proposition 1**.

### 2.3 Controlled supOU process

We consider a control problem of discharge by an environmental manager to mitigate an environmental and/or resource supply issue. For brevity, the dynamics are considered from $t=0$ instead of starting from a stationary state. The control $u=\{u_t\}_{>0}$ is assumed to be progressively measurable as usual, and assume the following controlled process $X$:

$$X_t=\underline{X}+\left(X_0-\underline{X}\right)\int_0^{+\infty}e^{-\lambda t}\pi(\mathrm{d}\lambda)+\int_0^t\int_0^{+\infty}e^{-\lambda(t-s)}\left(u_s\pi(\mathrm{d}\lambda)\mathrm{d}s+\int_0^{+\infty}z\mu(\mathrm{d}z,\mathrm{d}\lambda,\mathrm{d}s)\right),\quad t>0,\tag{21}$$

where the second term on the right-hand side of (21) appears due to the initial condition $X_0\geq 0$, and the second term includes the drift term to be controlled. The jump term is not controlled because it represents uncontrollable upstream inflow events.

The following assumptions have been made so that the objective functional is well-defined:

$$\limsup_{k\to+\infty}\frac{1}{kP}\mathbb{E}\left[\int_0^{kP}\left(X_s\right)^2\mathrm{d}s\right],\ \limsup_{k\to+\infty}\frac{1}{kP}\mathbb{E}\left[\int_0^{kP}\frac{1}{2}\left(X_s-\hat{X}_s\right)^2\mathrm{d}s\right],\ \limsup_{k\to+\infty}\frac{1}{kP}\mathbb{E}\left[\int_0^{kP}\frac{1}{2}u_s^2\mathrm{d}s\right]<+\infty\,.\tag{22}$$

Admissible set of controls, denoted as $\mathbb{U}$, is set as

$$\mathbb{U}=\left\{u=\left(u_t\right)_{t>0}\left|\begin{array}{c}u_t\in\mathbb{R}\text{ and progressibelly measurable at each }t>0\text{, satisfies (22),}\\ \text{and the controlled }X\text{ has a unique periodic proabbility measure.}\end{array}\right.\right\},\tag{23}$$

The assumption on the unique existence of the periodic probability measure is a technical one (See, also **Section 3.3** for the finite-dimensional case). Note that $\mathbb{U}$ is not null since $u\equiv 0$ belongs to it. Unless otherwise specified, we consider controls in $\mathbb{U}$.

Our objective functional is the LQ one:

$$J(u)=\limsup_{k\to+\infty}\frac{1}{kP}\mathbb{E}\left[\int_0^{kP}\left(\frac{1}{2}\left(X_s-\hat{X}_s\right)^2+\frac{w}{2}u_s^2\right)\mathrm{d}s\right]\tag{24}$$

with a period $P$ (e.g., one year) and $k\in\mathbb{N}$. In (24), the first term evaluates the deviation from some smooth and bounded target $\hat{X}_s\left(>\underline{X}\right)$ having the period $P$, while the second term evaluates the cost of control with $w>0$, a weighting coefficient balancing the two terms. This $J$ shall be minimized with respect to $u$. For example, it is a metric to be optimized for water resources management considering the seasonally varying environment and the river ecosystem. By construction, the discharge may become negative because we do not impose any state constraints. However, the occurrence probability of such events can be small if we set a sufficiently high target value or a strong penalization of control. The minimized $J$, an effective Hamiltonian, is defined as

$$H=\inf_u J(u).\tag{25}$$

We also consider a finite-dimensional counterpart for later use:

$$X_{n,t} = \underline{X} + \left(X_0 - \underline{X}\right)\sum_{i=1}^{n} c_i e^{-\lambda_i t} + \sum_{i=1}^{n}\int_0^t e^{-\lambda_i (t-s)}\left(c_i u_s \mathrm{d}s + \mathrm{d}L_s^{(i)}\right) = \underline{X} + \sum_{i=1}^{n} X_t^{(i)}, \quad t > 0 \tag{26}$$

with $X_t^{(i)}$ governed by

$$\mathrm{d}X_t^{(i)} = \left(-\lambda_i X_t^{(i)} + c_i u_t\right)\mathrm{d}t + \mathrm{d}L_t^{(i)}, \quad t > 0, \quad X_t^{(0)} = c_i X_0 . \tag{27}$$

The second term in the right-hand side of (26) is actually not so important as we consider a control problem in an infinite horizon. The objective functional $J_n$ and its optimized value $H_n$ are set as

$$J_n(u) = \limsup_{k \to +\infty} \frac{1}{kP}\mathbb{E}\left[\int_0^{kP}\left(\frac{1}{2}\left(X_{n,s} - \hat{X}_s\right)^2 + \frac{w}{2}u_s^2\right)\mathrm{d}s\right] \tag{28}$$

and

$$H_n = \inf_u J_n(u), \tag{29}$$

respectively. The admissible set of $u$ should be appropriately modified for the finite-dimensional case such that (22) with $X_s$ replaced by $X_{n,s}$ is satisfied.

**Remark 1:** Our formulation applies to a time-independent case where all model parameters are constants.

**Remark 2:** There is no technical difficulty to multiply $\left(X_{n,s} - \hat{X}_s\right)^2$ by a deterministic, smooth time-dependent, and positive weighting factor $w'_t$. We do not use this formulation except for **Section 4** for simplicity.

**Remark 3:** Our formulation assumes that the same control is applied to all $X^{(i)}$, which is a kind of finite-dimensional control against the infinite-dimensional system that is of importance in engineering optimization (Burger et al., 2020; Katz and Fridman, 2021). Instead, one may formulate a more complex control where different controls are applied to different $X_t^{(i)}$:

$$\mathrm{d}X_t^{(i)} = \left(-\lambda_i X_t^{(i)} + c_i u_{i,t}\right)\mathrm{d}t + \mathrm{d}L_t^{(i)}, \quad t \in \mathbb{R} . \tag{30}$$

The resulting HJB equation is different from ours but is still similar to it. The higher dimensional control would theoretically perform better, although it is more difficult to implement. This case is interesting by itself and will be analyzed elsewhere.

## 3 Control problem

### 3.1 Finite-dimensional case

We firstly discuss a finite-dimensional case as a first step. Then, as a second step, we consider an infinite-dimensional case as a formal limit of the finite-dimensional one. Their linkages are discussed as well.

### 3.1.1 Derivation of the Riccati equation

Based on a dynamic programming argument (e.g., Øksendal and Sulem, 2019; Yoshioka and Tsujimura, 2022), the HJB equation associated to (29) is derived as

$$-H_n+\frac{\partial\Phi}{\partial s}+\inf_{u\in\mathbb{R}}\left\{\begin{array}{l}\sum_{i=1}^{n}(-\lambda_i x_i+c_i u)\frac{\partial\Phi}{\partial x_i}+\frac{1}{2}\left(\sum_{i=1}^{n}x_i+\underline{X}-\hat{X}_s\right)^2\\+\frac{w}{2}u^2+\sum_{i=1}^{n}c_i\int_0^{+\infty}\Delta_i\Phi(z)v(\mathrm{d}z)\end{array}\right\}=0 \quad (31)$$

with

$$\Delta_i\Phi(z)=\Phi(s,x_1,x_2,...,x_i+z,...,x_n)-\Phi(s,x_1,x_2,...,x_n). \quad (32)$$

Here, $\Phi:[0,P]\times\mathbb{R}_+^n\to\mathbb{R}$ is a sufficiently smooth potential function having the period $P$ with respect to the first argument, and $H_n\in\mathbb{R}$ is identified as an effective Hamiltonian of (29) with abusing notations. A solution to (31) is a couple $(H_n,\Phi)$.

The HJB equation (31) is rearranged as

$$-H_n+\sum_{i=1}^{n}(-\lambda_i x_i)\frac{\partial\Phi}{\partial x_i}-\frac{1}{2w}\left(\sum_{i=1}^{n}c_i\frac{\partial\Phi}{\partial x_i}\right)^2+\sum_{i=1}^{n}c_i\int_0^{+\infty}\Delta_i\Phi(z)v(\mathrm{d}z)+\frac{1}{2}\left(\sum_{i=1}^{n}x_i+\underline{X}-\hat{X}_s\right)^2=0 \quad (33)$$

with (a candidate of) the optimal control

$$u^*(s,x)=\arg\min_{\mathbb{R}}\left\{\sum_{i=1}^{n}c_i u\frac{\partial\Phi}{\partial x_i}+\frac{w}{2}u_s^2\right\}=-\frac{1}{w}\sum_{i=1}^{n}c_i\frac{\partial\Phi}{\partial x_i}. \quad (34)$$

We guess a quadratic $\Phi$ with time-dependent coefficients, which are a symmetric matrix $\{A_{i,j}\}_{1\le i,j\le n}$, a vector $\{B_i\}_{1\le i\le n}$, and a scalar $E$:

$$\Phi=\frac{1}{2}\sum_{i,j=1}^{n}A_{i,j}x_i x_j+\sum_{i=1}^{n}B_i x_i+E. \quad (35)$$

Substituting (35) into each term of (33) yields the following identities:

$$\sum_{i=1}^{n}(-\lambda_i x_i)\frac{\partial\Phi}{\partial x_i}=-\sum_{i,j=1}^{n}\frac{1}{2}(\lambda_i+\lambda_j)A_{i,j}x_i x_j-\sum_{i=1}^{n}\lambda_i B_i x_i, \quad (36)$$

$$\begin{aligned}-\frac{1}{2w}\left(\sum_{i=1}^{n}c_i\frac{\partial\Phi}{\partial x_i}\right)^2&=-\frac{1}{2w}\left(\sum_{i=1}^{n}c_i\left(\sum_{j=1}^{n}A_{i,j}x_j+B_i\right)\right)^2\\&=-\frac{1}{2w}\sum_{i,j=1}^{n}D_iD_jx_ix_j-\frac{1}{w}\left(\sum_{j=1}^{n}c_jB_j\right)\sum_{i=1}^{n}D_ix_i-\frac{1}{2w}\left(\sum_{i=1}^{n}c_iB_i\right)^2\end{aligned} \quad (37)$$

with $D_i=\sum_{j=1}^{n}c_jA_{i,j}$,

$$\sum_{i=1}^{n}c_i\int_0^{+\infty}\Delta_i\Phi(z)v(\mathrm{d}z)=M_1\sum_{i=1}^{n}D_ix_i+\frac{1}{2}M_2\sum_{i=1}^{n}c_iA_{i,i}+M_1\sum_{i=1}^{n}c_iB_i \quad (38)$$

by

$$\begin{aligned}\Delta_i\Phi(z) &= \frac{1}{2}\sum_{k,l=1}^{n} A_{k,l}\left(x_k+\delta_{i,k}z\right)\left(x_l+\delta_{i,l}z\right)+\sum_{k=1}^{n} B_k\left(x_k+\delta_{i,k}z\right)-\left(\frac{1}{2}\sum_{k,l=1}^{n} A_{k,l}x_k x_l+\sum_{k=1}^{n} B_k x_k\right)\\ &= \frac{1}{2}\sum_{k,l=1}^{n} A_{k,l}\left(\delta_{i,k}zx_l+\delta_{i,l}zx_k+\delta_{i,k}\delta_{i,l}z^2\right)+\sum_{k=1}^{n} B_k\delta_{i,k}z\\ &= \frac{1}{2}A_{i,i}z^2+\sum_{k=1}^{n} A_{k,i}x_k z+B_i z\end{aligned} \tag{39}$$

with the Kronecker Delta $\delta_{i,j}$, and

$$\frac{1}{2}\left(\sum_{i=1}^{n} x_i+\underline{X}-\hat{X}_s\right)^2=\frac{1}{2}\left(\sum_{i=1}^{n} x_i-\bar{X}_s\right)^2=\frac{1}{2}\sum_{i,j=1}^{n} x_i x_j-\bar{X}_s\sum_{i=1}^{n} x_i+\frac{1}{2}\left(\bar{X}_s\right)^2 \tag{40}$$

with $\bar{X}_s=\hat{X}_s-\underline{X}>0$. Consequently, (33) becomes

$$\begin{aligned}0 &= \sum_{i,j=1}^{n}\left(\frac{\partial A_{i,j}}{\partial s}-\frac{1}{2}\left(\lambda_i+\lambda_j\right)A_{i,j}-\frac{1}{2w}D_i D_j+\frac{1}{2}\right)x_i x_j\\ &\quad+\sum_{i=1}^{n}\left(\frac{\partial B_i}{\partial s}-\lambda_i B_i-\frac{1}{w}D_i\sum_{j=1}^{n} c_j B_j+M_1 D_i-\bar{X}_s\right)x_i\\ &\quad-H_n+\frac{\partial E}{\partial s}+\frac{1}{2}M_2\sum_{i=1}^{n} c_i A_{i,i}+M_1\sum_{i=1}^{n} c_i B_i-\frac{1}{2w}\left(\sum_{i=1}^{n} c_i B_i\right)^2+\frac{1}{2}\left(\bar{X}_s\right)^2\end{aligned} \tag{41}$$

We then obtain the time-backward finite-dimensional Riccati equation governing $A_{i,j}$, $B_i$, $E$, $H_n$:

$$0=\frac{\partial A_{i,j}}{\partial s}-\left(\lambda_i+\lambda_j\right)A_{i,j}-\frac{1}{w}D_i D_j+1 \quad (1\le i,j\le n,\ 0\le s\le P), \tag{42}$$

$$0=\frac{\partial B_i}{\partial s}-\lambda_i B_i-\frac{D_i}{w}\sum_{j=1}^{n} c_j B_j+M_1 D_i-\bar{X}_s \quad (1\le i\le n,\ 0\le s\le P), \tag{43}$$

$$0=\frac{\partial E}{\partial s}-H_n-\frac{1}{2w}\left(\sum_{i=1}^{n} c_i B_i\right)^2+\frac{1}{2}M_2\sum_{i=1}^{n} c_i A_{i,i}+M_1\sum_{i=1}^{n} c_i B_i+\frac{1}{2}\left(\bar{X}_s\right)^2 \quad (0\le s\le P). \tag{44}$$

(A candidate of) The effective Hamiltonian $H_n$ is then obtained as

$$H_n=\frac{1}{P}\int_0^P\left\{-\frac{1}{2w}\left(\sum_{i=1}^{n} c_i B_i\right)^2+\frac{1}{2}M_2\sum_{i=1}^{n} c_i A_{i,i}+M_1\sum_{i=1}^{n} c_i B_i+\frac{1}{2}\left(\bar{X}_s\right)^2\right\}\mathrm{d}s. \tag{45}$$

The total numbers of equations (42), (43), (44) (or (45)) and the unknowns equal $n^2+n+1$. The equations (42)–(44) are integrated backward in time. Finally, (34) becomes

$$u^*=-\frac{1}{w}\sum_{i=1}^{n} c_i\frac{\partial\Phi}{\partial x_i}=-\frac{1}{w}\sum_{i=1}^{n} c_i\left(\sum_{j=1}^{n} A_{i,j}x_j+B_i\right). \tag{46}$$

Consequently, we reduced a multi-dimensional HJB equation to a finite-dimensional Riccati equation.

### 3.1.2 Existence and uniqueness of the Riccati equation

Our Riccati equation is the assembly of (42), (43), and (44) (or (45)), where equation (42) is a nonlinear system of $A$, equation (43) is a linear system of $B$ that gives $A$, and equation (44) is a trivial equation that gives $A, B$. We can solve the Riccati equation in this order.

The unique existence of a time-periodic positive semi-definite $A$ solving (42) is obtained from the argument similar to Theorem 6 of Ha and Choi (2019). Our noise processes differ from those used in their work, while the Riccati equation for the matrix parts is similar. The proof of **Proposition 2** below stating the unique solvability is provided in **Appendix C.2**.

***Proposition 2***

*There system consisting of (42), (43), (45) is uniquely solved by a smooth and positive semi-definite $n \times n$-dimensional matrix $A$ having the period $P$, $n$-dimensional smooth vector $B_i$ having the period $P$, and a constant $H_n$.*

### 3.2 Derivation of the Riccati equation for the infinite-dimensional case

An infinite-dimensional limit of the HJB equation is heuristically obtained. Its solution is defined as a functional of elements in a Hilbert space. Set the weighted Lebesgue space $L^p(\pi)$ ($p \geq 1$) as a collection of $\phi:(0,+\infty) \to \mathbb{R}$ equipped with the norm

$$\|\phi\|_{L^p(\pi)} = \left(\int_0^{+\infty} |\phi^p(\lambda)| \pi(\mathrm{d}\lambda)\right)^{\frac{1}{p}} < +\infty \quad \text{if} \quad p < +\infty, \tag{47}$$

$$\|\phi\|_{L^\infty(\pi)} = \sup_{\lambda > 0} \phi(\lambda) \quad \text{if} \quad p = +\infty. \tag{48}$$

The constant function $\phi(\lambda) \equiv 1$ belongs to $L^p(\pi)$ ($p \geq 1$). Furthermore $\phi \in L^1(\pi)$, if $\phi \in L^2(\pi)$:

$$\|\phi\|_{L^1(\pi)} = \int_0^{+\infty} |\phi(\lambda)| \pi(\mathrm{d}\lambda) \leq \sqrt{\int_0^{+\infty} (\phi(\lambda))^2 \pi(\mathrm{d}\lambda)} \sqrt{\int_0^{+\infty} 1 \pi(\mathrm{d}\lambda)} \leq \|\phi\|_{L^2(\pi)}. \tag{49}$$

In the split of a Markovian lift, we rewrite the controlled supOU process as

$$X_t = \underline{X} + \int_0^{+\infty} Y_t(\lambda) \pi(\mathrm{d}\lambda), \quad t > 0 \tag{50}$$

with the auxiliary stochastic field

$$Y_t(\lambda) = (X_0 - \underline{X}) e^{-\lambda t} + \int_0^t e^{-\lambda(t-s)} \left(u_s \mathrm{d}s + \mathrm{d}L_s(\lambda)\right), \quad t \in \mathbb{R}, \quad \lambda > 0, \tag{51}$$

where $(L_t(\lambda))_{t \in \mathbb{R}}$ is a space–time Lévy process (e.g., Dunst et al., 2012; Barth and Stein, 2018) such that formally, $\mathbb{E}[\mathrm{d}L_s(\cdot)] = M_1 \mathrm{d}s$ and $\mathbb{E}[\mathrm{d}L_s(\lambda) \mathrm{d}L_s(\theta)] = \left(\chi_{\{\theta \neq \lambda\}} M_1^2 (\mathrm{d}s)^2 + \hat{\delta}_{\{\theta = \lambda\}} \chi_{\{\theta = \lambda\}} M_2 \mathrm{d}s\right)$ with $\hat{\delta}_{\{\omega = (\cdot)\}} = \delta_{\{\omega = (\cdot)\}} \dfrac{\mathrm{d}\omega}{\pi(\mathrm{d}\omega)}$, where $\delta$ is the Dirac Delta. In this way, the supOU process is based on a stochastic partial differential equation.

**Remark 4:** Stochastic partial differential equations driven by space–time white noises are not always well-defined (e.g., Kuehn and Riedler (2014)), while the supOU process avoids this issue if the conditions (2)–(3) are satisfied. They guarantee that the density $\pi$ decays sufficiently fast as $\lambda \to +\infty$ and is sufficiently regular near $\lambda = 0$. The correlation of $\mathrm{d}L$ is consistent with the HJB equation (**Appendix C.3**).

Using (50)–(51), the objective functional becomes

$$J(u) = \limsup_{k\to+\infty} \frac{1}{kP} \mathbb{E}\left[ \int_0^{kP} \left( \frac{1}{2} \left( \int_0^{+\infty} Y_{u,s}(\lambda)\pi(\mathrm{d}\lambda) - \bar{X}_s \right)^2 + \frac{w}{2} u_s^2 \right) \mathrm{d}s \right]. \tag{52}$$

The requirement of the square integrability of $Y$ comes from the quadratic term in (52). The potential is denoted as $V = V(s, Y)$, which is a function of $Y \in L^2(\pi)$ and $s \in [0, P]$. Its dependence on $s$ is sometimes suppressed below for simplicity. Given $Y \in L^2(\pi)$, each term of (31) is expected to converge as $n \to +\infty$:

$$\sum_{i=1}^n (-\lambda_i x_i + c_i u) \frac{\partial \Phi}{\partial x_i} + \frac{w}{2} u^2 \to \int_0^{+\infty} (-\lambda Y(\lambda) + u) \nabla V(Y) \pi(\mathrm{d}\lambda) + \frac{w}{2} u^2, \tag{53}$$

$$\frac{1}{2}\left( \sum_{i=1}^n x_i + \underline{X} - \hat{X} \right)^2 \to \frac{1}{2} \int_0^{+\infty}\int_0^{+\infty} Y(\lambda)Y(\theta)\pi(\mathrm{d}\lambda)\pi(\mathrm{d}\theta) - \bar{X}\int_0^{+\infty} Y(\lambda)\pi(\mathrm{d}\lambda) + \frac{1}{2}\bar{X}^2, \tag{54}$$

$$\sum_{i=1}^n c_i \int_0^{+\infty} \Delta_i \Phi(z) v(\mathrm{d}z) \to \int_0^{+\infty} \left( \int_0^{+\infty} V\left( Y(\cdot) + z\hat{\delta}_{\{\omega=(\cdot)\}} \right) \pi(\mathrm{d}\omega) - V(Y) \right) v(\mathrm{d}z), \tag{55}$$

and $H_n \to H$. Here, $\nabla V$ is the Fréchet derivative of $V$ (Definition 2.9 of Serovajsky, 2017). We then expect an optimal control

$$u^*(s, Y) = \arg\min_{u\in\mathbb{R}} \left\{ u \int_0^{+\infty} \nabla V(Y)\pi(\mathrm{d}\lambda) + \frac{w}{2}u^2 \right\} = -\frac{1}{w}\int_0^{+\infty} \nabla V(Y)\pi(\mathrm{d}\lambda) \tag{56}$$

and the infinite-dimensional HJB equation

$$\begin{aligned} &-H + \frac{\partial V}{\partial s} - \int_0^{+\infty} \lambda \nabla V(Y) Y(\lambda)\pi(\mathrm{d}\lambda) - \frac{1}{2w}\left( \int_0^{+\infty} \nabla V(Y)\pi(\mathrm{d}\lambda) \right)^2 \\ &+ \int_0^{+\infty} \left( \int_0^{+\infty} V\left( Y(\cdot) + z\hat{\delta}_{\{\omega=(\cdot)\}} \right)\pi(\mathrm{d}\omega) - V(Y) \right) v(\mathrm{d}z) \\ &+ \frac{1}{2}\int_0^{+\infty}\int_0^{+\infty} Y(\lambda)Y(\theta)\pi(\mathrm{d}\lambda)\pi(\mathrm{d}\theta) - \bar{X}\int_0^{+\infty} Y(\lambda)\pi(\mathrm{d}\lambda) + \frac{1}{2}\bar{X}^2 = 0 \end{aligned} . \tag{57}$$

We guess $V$ as

$$\begin{aligned} V(s, Y) = &\frac{1}{2}\int_0^{+\infty}\int_0^{+\infty} \Gamma_s(\theta, \lambda) Y(\lambda) Y(\theta)\pi(\mathrm{d}\lambda)\pi(\mathrm{d}\theta) \\ &+ \int_0^{+\infty} \gamma_s(\lambda) Y(\lambda)\pi(\mathrm{d}\lambda) + R_s \end{aligned}, \quad Y \in L^2(\pi), \quad s \in [0, P] \tag{58}$$

with a time-dependent positive-semidefinite $\pi \otimes \pi$-integrable $\Gamma$ and $\pi$-integrable $\gamma$; we say that $\Gamma : (0, +\infty)^2 \to \mathbb{R}$ belongs to $L^1(\pi \otimes \pi)$ if

$$\|\Gamma\|_{L^1(\pi\otimes\pi)} = \int_0^{+\infty}\int_0^{+\infty} |\Gamma(\theta, \lambda)| \pi(\mathrm{d}\lambda)\pi(\mathrm{d}\theta) < +\infty. \tag{59}$$

We say that $\Gamma:(0,+\infty)^2 \to \mathbb{R}$ (resp., $\gamma:(0,+\infty)\to\mathbb{R}$) belongs to $L^\infty(\pi\otimes\pi)$ (resp., $L^\infty(\pi)$) if it is bounded a.e. on $\pi$.

For $\Gamma\in L^\infty(\pi\otimes\pi)$ and $\gamma\in L^\infty(\pi)$, we obtain

$$\begin{aligned}\nabla V(Y)Z &= \lim_{h\to 0}\frac{V(Y+hZ)-V(Y)}{h} \\ &= \int_0^{+\infty}\left(\int_0^{+\infty}\Gamma(\theta,\lambda)Y(\theta)\pi(\mathrm{d}\theta)+\gamma(\lambda)\right)Z(\lambda)\pi(\mathrm{d}\lambda)\end{aligned}, \quad Y,Z\in L^2(\pi). \tag{60}$$

Namely,

$$\nabla V(Y)=\int_0^{+\infty}\Gamma(\theta,\lambda)Y(\theta)\pi(\mathrm{d}\theta)+\gamma(\lambda),\quad Y\in L^2(\pi), \tag{61}$$

where we used the symmetry of $\Gamma$. We obtain

$$\begin{aligned}&-\int_0^{+\infty}\lambda\nabla V(Y)Y(\lambda)\pi(\mathrm{d}\lambda)\\ &=-\frac{1}{2}\int_0^{+\infty}\int_0^{+\infty}(\theta+\lambda)\Gamma(\theta,\lambda)Y(\theta)Y(\lambda)\pi(\mathrm{d}\theta)\pi(\mathrm{d}\lambda)-\int_0^{+\infty}\lambda\gamma(\lambda)Y(\lambda)\pi(\mathrm{d}\lambda)\end{aligned}, \tag{62}$$

$$\begin{aligned}&-\frac{1}{2w}\left(\int_0^{+\infty}\nabla V(Y)\pi(\mathrm{d}\lambda)\right)^2\\ &=-\frac{1}{2w}\left(\begin{array}{l}\int_0^{+\infty}\int_0^{+\infty}\left(\int_0^{+\infty}\Gamma(\theta_1,\lambda_1)\pi(\mathrm{d}\theta_1)\int_0^{+\infty}\left(\Gamma(\theta_2,\lambda_2)\pi(\mathrm{d}\theta_2)\right)\right)Y(\lambda_1)Y(\lambda_2)\pi(\mathrm{d}\lambda_1)\pi(\mathrm{d}\lambda_2)\\ +\left(\int_0^{+\infty}\gamma(\lambda)\pi(\mathrm{d}\lambda)\right)^2+2\left(\int_0^{+\infty}\gamma(\lambda)\pi(\mathrm{d}\lambda)\right)\left(\int_0^{+\infty}\left(\int_0^{+\infty}\Gamma(\theta,\lambda)\pi(\mathrm{d}\theta)\right)Y(\lambda)\pi(\mathrm{d}\lambda)\right)\end{array}\right),\end{aligned} \tag{63}$$

and

$$\begin{aligned}&\int_0^{+\infty}V\left(Y(\cdot)+z\hat{\delta}_{\{\omega=(\cdot)\}}\right)\pi(\mathrm{d}\omega)-V(Y)\\ &=\frac{1}{2}\int_0^{+\infty}\int_0^{+\infty}\int_0^{+\infty}\Gamma(\theta,\lambda)\left(Y(\lambda)+z\hat{\delta}_{\{\omega=\lambda\}}\right)\left(Y(\theta)+z\hat{\delta}_{\{\omega=\theta\}}\right)\pi(\mathrm{d}\lambda)\pi(\mathrm{d}\theta)\pi(\mathrm{d}\omega)\\ &+\int_0^{+\infty}\int_0^{+\infty}\gamma(\lambda)\left(Y(\lambda)+z\hat{\delta}_{\{\omega=\lambda\}}\right)\pi(\mathrm{d}\lambda)\pi(\mathrm{d}\omega)+R\\ &-\left(\frac{1}{2}\int_0^{+\infty}\int_0^{+\infty}\Gamma(\theta,\lambda)Y(\lambda)Y(\theta)\pi(\mathrm{d}\lambda)\pi(\mathrm{d}\theta)+\int_0^{+\infty}\gamma(\lambda)Y(\lambda)\pi(\mathrm{d}\lambda)+R\right)\\ &=\frac{1}{2}z^2\int_0^{+\infty}\int_0^{+\infty}\int_0^{+\infty}\Gamma(\theta,\lambda)\hat{\delta}_{\{\omega=\lambda\}}\hat{\delta}_{\{\omega=\theta\}}\pi(\mathrm{d}\lambda)\pi(\mathrm{d}\theta)\pi(\mathrm{d}\omega)\\ &+\frac{1}{2}z\int_0^{+\infty}\int_0^{+\infty}\int_0^{+\infty}\Gamma(\theta,\lambda)\left(\hat{\delta}_{\{\omega=\theta\}}Y(\lambda)+\hat{\delta}_{\{\omega=\lambda\}}Y(\theta)\right)\pi(\mathrm{d}\lambda)\pi(\mathrm{d}\theta)\pi(\mathrm{d}\omega)+z\int_0^{+\infty}\gamma(\lambda)\pi(\mathrm{d}\lambda)\\ &=\frac{1}{2}z^2\int_0^{+\infty}\Gamma(\theta,\theta)\pi(\mathrm{d}\theta)+z\int_0^{+\infty}\left(\int_0^{+\infty}\Gamma(\theta,\lambda)\pi(\mathrm{d}\theta)\right)Y(\lambda)\pi(\mathrm{d}\lambda)+z\int_0^{+\infty}\gamma(\lambda)\pi(\mathrm{d}\lambda)\end{aligned}, \tag{64}$$

and thus

$$\begin{aligned}&\int_0^{+\infty}\left(\int_0^{+\infty}V\left(Y(\cdot)+z\hat{\delta}_{\{\omega=(\cdot)\}}\right)\pi(\mathrm{d}\lambda)\pi(\mathrm{d}\omega)-V(Y)\right)v(\mathrm{d}z)\\ &=M_1\int_0^{+\infty}\left(\int_0^{+\infty}\Gamma(\theta,\lambda)\pi(\mathrm{d}\theta)\right)Y(\lambda)\pi(\mathrm{d}\lambda)+M_1\int_0^{+\infty}\gamma(\lambda)\pi(\mathrm{d}\lambda)+\frac{1}{2}M_2\int_0^{+\infty}\Gamma(\theta,\theta)\pi(\mathrm{d}\theta)\end{aligned}. \tag{65}$$

In summary, (57) reduces to

$$
\begin{aligned}
0 &= -H + \frac{\partial R}{\partial s} + M_1\int_0^{+\infty}\gamma(\lambda)\pi(\mathrm{d}\lambda) + \frac{1}{2}M_2\int_0^{+\infty}\Gamma(\theta,\theta)\pi(\mathrm{d}\theta) - \frac{1}{2w}\left(\int_0^{+\infty}\gamma(\lambda)\pi(\mathrm{d}\lambda)\right)^2 + \frac{1}{2}\bar{X} \\
&+ \int_0^{+\infty}\left(\begin{array}{l}\frac{\partial \gamma}{\partial s} - \lambda\gamma(\lambda) + M_1\int_0^{+\infty}\Gamma(\theta,\lambda)\pi(\mathrm{d}\theta) \\ -\frac{1}{w}\int_0^{+\infty}\gamma(\theta)\pi(\mathrm{d}\theta)\int_0^{+\infty}\Gamma(\theta,\lambda)\pi(\mathrm{d}\theta) - \bar{X}\end{array}\right)Y(\lambda)\pi(\mathrm{d}\lambda) \\
&+ \int_0^{+\infty}\int_0^{+\infty}\left\{\begin{array}{l}\frac{1}{2}\frac{\partial \Gamma}{\partial s} - \frac{1}{2}(\theta+\lambda)\Gamma(\theta,\lambda) \\ -\frac{1}{2w}\int_0^{+\infty}\Gamma(\theta_1,\lambda)\pi(\mathrm{d}\theta_1)\int_0^{+\infty}\Gamma(\theta_2,\theta)\pi(\mathrm{d}\theta_2) + \frac{1}{2}\end{array}\right\}Y(\theta)Y(\lambda)\lambda\pi(\mathrm{d}\theta)\pi(\mathrm{d}\lambda)
\end{aligned}
\quad . \tag{66}
$$

From (66), we obtain the identities of the integral operator Riccati equation solved backward in time:

$$
\begin{aligned}
0 &= \frac{\partial \Gamma(\theta,\lambda)}{\partial s} - (\theta+\lambda)\Gamma(\theta,\lambda) \\
&\quad - \frac{1}{w}\int_0^{+\infty}\Gamma(\theta_1,\lambda)\pi(\mathrm{d}\theta_1)\int_0^{+\infty}\Gamma(\theta_2,\theta)\pi(\mathrm{d}\theta_2) + 1
\end{aligned}
\qquad (\lambda,\theta > 0,\ s\in[0,P]), \tag{67}
$$

$$
\begin{aligned}
0 &= \frac{\partial \gamma(\lambda)}{\partial s} - \lambda\gamma(\lambda) \\
&\quad - \frac{1}{w}\int_0^{+\infty}\Gamma(\theta,\lambda)\pi(\mathrm{d}\theta)\int_0^{+\infty}\gamma(\lambda)\pi(\mathrm{d}\lambda) + M_1\int_0^{+\infty}\Gamma(\theta,\lambda)\pi(\mathrm{d}\theta) - \bar{X}
\end{aligned}
\qquad (\lambda > 0,\ s\in[0,P]), \tag{68}
$$

$$
\begin{aligned}
0 &= \frac{\partial R}{\partial s} - H - \frac{1}{2w}\left(\int_0^{+\infty}\gamma(\lambda)\pi(\mathrm{d}\lambda)\right)^2 \\
&\quad + \frac{1}{2}M_2\int_0^{+\infty}\Gamma(\theta,\theta)\pi(\mathrm{d}\theta) + M_1\int_0^{+\infty}\gamma(\lambda)\pi(\mathrm{d}\lambda) + \frac{1}{2}\bar{X}^2
\end{aligned}
\qquad (s\in[0,P]). \tag{69}
$$

From (69), we obtain the effective Hamiltonian

$$
H = \frac{1}{P}\int_0^P\left(\frac{1}{2}M_2\int_0^{+\infty}\Gamma(\theta,\theta)\pi(\mathrm{d}\theta) + M_1\int_0^{+\infty}\gamma(\lambda)\pi(\mathrm{d}\lambda) - \frac{1}{2w}\left(\int_0^{+\infty}\gamma(\lambda)\pi(\mathrm{d}\lambda)\right)^2 + \frac{1}{2}\bar{X}^2\right)\mathrm{d}s\,. \tag{70}
$$

The optimal control $u^* = u_s^*$ then would become

$$
u_s^* = -\frac{1}{w}\int_0^{+\infty}\left(\int_0^{+\infty}\Gamma(\theta,\lambda)Y(\theta)\pi(\mathrm{d}\theta) + \gamma(\lambda)\right)\pi(\mathrm{d}\lambda),\ \ Y\in L^2(\pi),\ \ s\in[0,P]. \tag{71}
$$

Consequently, we reduced an infinite-dimensional HJB equation to a Riccati equation of the integro-partial differential form.

**Remark 5**: Each term in (67)–(69) has the corresponding discrete analog in (42)–(44), suggesting their consistency.

### 3.3 On the stability of the controlled process

A concern of controlled dynamics is their dissipativity and the unique existence of a periodic measure. For the finite-dimensional case, set $|\mathbf{x}| = \sqrt{\sum_{i=1}^{n} x_i^2}$ and

$$d_i(\mathbf{x}) = -\lambda_i x_i - \frac{1}{w} c_i \sum_{j,k=1}^{n} c_k A_{k,j} x_j \quad (1 \le i \le n),\ \ \mathbf{x} \in \mathbb{R}^n \ . \tag{72}$$

We then have

$$\mathbf{x} \cdot d_i(\mathbf{x}) = -\sum_{i=1}^{n} \lambda_i x_i^2 - \frac{1}{w} \sum_{i=1}^{n} c_i x_i \sum_{j,k=1}^{n} c_k A_{k,j} x_j \le -\sum_{i=1}^{n} \lambda_i x_i^2 + \frac{1}{w} \sqrt{\sum_{i=1}^{n} c_i^2} \sqrt{\sum_{j=1}^{n} \sum_{k=1}^{n} \left(c_k A_{k,j}\right)^2} |\mathbf{x}|^2 . \tag{73}$$

The controlled system is an OU-type; i.e., it has a unique path-wise solution global in time. By Example 4.1 of Guo and Sun (2021), if $\frac{1}{w} \sqrt{\sum_{i=1}^{n} c_i^2} \sqrt{\sum_{j=1}^{n} \sum_{k=1}^{n} \left(c_k A_{k,j}\right)^2}$ is small (i.e., control cost is highly penalized), then $C > 0$ and the controlled dynamics are dissipative, guaranteeing the sufficient condition of Theorem 3.14 of the same literature and hence the unique existence of a periodic probability measure.

For the infinite-dimensional case, the situation is different. We have

$$\begin{aligned}
&\int_0^{+\infty} Y(\lambda)\left(-\lambda Y(\lambda) - \frac{1}{w}\int_0^{+\infty}\int_0^{+\infty} \Gamma(\theta,\omega) Y(\theta)\pi(\mathrm{d}\theta)\pi(\mathrm{d}\omega)\right)\pi(\mathrm{d}\lambda) \\
&\le -\int_0^{+\infty} \lambda \left(Y(\lambda)\right)^2 \pi(\mathrm{d}\lambda) \\
&+ \frac{1}{w}\int_0^{+\infty} |Y(\lambda)| \left(\int_0^{+\infty}\left(\int_0^{+\infty} \Gamma(\theta,\omega)\pi(\mathrm{d}\omega)\right)|Y(\theta)|\pi(\mathrm{d}\theta)\right)\pi(\mathrm{d}\lambda) \\
&\le -\int_0^{+\infty} \lambda \left(Y(\lambda)\right)^2 \pi(\mathrm{d}\lambda) \\
&+ \frac{1}{w}\int_0^{+\infty} |Y(\lambda)| \sqrt{\int_0^{+\infty}\left(\int_0^{+\infty} |\Gamma(\theta,\omega)|\pi(\mathrm{d}\omega)\right)^2 \pi(\mathrm{d}\theta)} \|Y\|_{L^2(\pi)} \pi(\mathrm{d}\lambda) \\
&\le -\int_0^{+\infty} \lambda \left(Y(\lambda)\right)^2 \pi(\mathrm{d}\lambda) + \frac{C_1}{w}\|Y\|_{L^2(\pi)}^2
\end{aligned} \quad ,\ Y \in L^2(\pi) \tag{74}$$

with

$$C_1 = \sqrt{\int_0^{+\infty}\left(\int_0^{+\infty} |\Gamma(\theta,\omega)|\pi(\mathrm{d}\omega)\right)^2 \pi(\mathrm{d}\theta)} \ge 0 . \tag{75}$$

If there is a constant $C' > 0$ such that

$$-\int_0^{+\infty} \lambda\left(Y(\lambda)\right)^2 \pi(\mathrm{d}\lambda) + \frac{C_1}{w}\|Y\|_{L^2(\pi)}^2 \le -C'\|Y\|_{L^2(\pi)}^2 \ \text{ or } \ \left(C' + \frac{C_1}{w}\right)\|Y\|_{L^2(\pi)}^2 \le \int_0^{+\infty} \lambda\left(Y(\lambda)\right)^2 \pi(\mathrm{d}\lambda), \tag{76}$$

then, the drift of the controlled dynamics is dissipative (e.g., Theorem 16.5 of Peszat and Zabczyk (2007)), suggesting their stochastic (exponential) stability.

Unfortunately, (76) is broken for $Y$ having a large spike at $0 < \lambda << C' + \frac{C_1}{w}$. This can be related to the sub-exponential autocorrelation (11) arising from the distributed reversion speed $\lambda \in (0, +\infty)$ with which the dynamics transit slower than an exponential speed. Therefore, we need a different metric other than dissipativity to evaluate the performance of the control, and we consider a KBE in **Section 3.4**.

### 3.4 Performance evaluation in a finite-dimensional case

Evaluating the performance of controls usually requires sampling a large number of controlled stochastic paths, which can be time-consuming. Our LQ framework allows for the derivation of a KBE for this purpose. This KBE is also useful in evaluating the performance of optimal controls, at least computationally, when the conventional analysis based on dissipativity fails.

We evaluate the performance of the optimal control $u^*$:

$$C = \limsup_{k\to+\infty} \frac{1}{kP} \mathbb{E}\left[ \int_0^{kP} \frac{1}{2}\left(u_s^*\right)^2 \mathrm{d}s \right]. \tag{77}$$

This is the controlling cost with $u = u^*$. The other term measuring the deviation of the state from the target can also be evaluated, but this term can be evaluated by finding $C$, as shown in **Section 4.5.3**. As in the beginning of **Section 3.1.1**, we deduce the KBE

$$\begin{aligned} &-C + \frac{\partial \Psi}{\partial s} + \sum_{i=1}^{n}\left( -\lambda_i x_i - c_i \frac{1}{w}\sum_{k=1}^{n} c_k \left( \sum_{l=1}^{n} A_{k,l} x_l + B_k \right)\right)\frac{\partial \Psi}{\partial x_i} + \sum_{i=1}^{n} c_i \int_0^{+\infty} \Delta_i \Psi(z) v(\mathrm{d}z) \\ &+ \frac{1}{2}\left( -\frac{1}{w}\sum_{i=1}^{n} c_i \left( \sum_{j=1}^{n} A_{i,j} x_j + B_i \right)\right)^2 = 0 \end{aligned}. \tag{78}$$

We guess the quadratic solution with time-dependent coefficients having the period $P$, which are a symmetric matrix $\left\{S_{i,j}\right\}_{1\le i,j\le n}$, a vector $\left\{N_i\right\}_{1\le i\le n}$, and a scalar $\tilde{E}$:

$$\Psi = \frac{1}{2}\sum_{i,j=1}^{n} S_{i,j} x_i x_j + \sum_{i=1}^{n} N_i x_i + \tilde{E}. \tag{79}$$

We then have

$$\sum_{i=1}^{n} \left(-\lambda_i x_i\right) \frac{\partial \Psi}{\partial x_i} = -\sum_{i,j=1}^{n} \frac{1}{2}\left(\lambda_i + \lambda_j\right) S_{i,j} x_i x_j - \sum_{i=1}^{n} \lambda_i N_i x_i, \tag{80}$$

$$\sum_{i=1}^{n} c_i \int_0^{+\infty} \Delta_i \Psi(z) v(\mathrm{d}z) = M_1 \sum_{i=1}^{n} F_i x_i + \frac{1}{2} M_2 \sum_{i=1}^{n} c_i S_{i,i} + M_1 \sum_{i=1}^{n} c_i N_i \tag{81}$$

with $F_i = \sum_{i=1}^{n} c_j S_{i,j}$,

$$\begin{aligned} &\sum_{i=1}^{n}\left( -\frac{1}{w} c_i \sum_{k=1}^{n} c_k \left( \sum_{l=1}^{n} A_{k,l} x_l + B_k \right)\right)\frac{\partial \Psi}{\partial x_i} \\ &= -\frac{1}{2w}\sum_{i,j=1}^{n}\left( D_i F_j + D_j F_i \right) x_i x_j \\ &-\frac{1}{w}\left( \sum_{k=1}^{n} c_k B_k \right)\sum_{i=1}^{n} F_i x_i - \frac{1}{w}\left( \sum_{i=1}^{n} c_i N_i \right)\sum_{i=1}^{n} D_i x_i - \frac{1}{w}\left( \sum_{k=1}^{n} c_k B_k \right)\left( \sum_{i=1}^{n} c_i N_i \right) \end{aligned}, \tag{82}$$

and

$$\frac{1}{2}\left( -\frac{1}{w}\sum_{i=1}^{n} c_i \left( \sum_{j=1}^{n} A_{i,j} x_j + B_i \right)\right)^2 = \frac{1}{2w^2}\sum_{i,j=1}^{n} D_i D_j x_i x_j + \frac{1}{w^2}\left( \sum_{i=1}^{n} c_i B_i \right)\sum_{i=1}^{n} D_i x_i + \frac{1}{2w^2}\left( \sum_{i=1}^{n} c_i B_i \right). \tag{83}$$

Consequently, we obtain

$$
\begin{aligned}
&\sum_{i,j=1}^{n}\left\{\frac{1}{2}\frac{\partial S_{i,j}}{\partial s}-\frac{1}{2}\left(\lambda_i+\lambda_j\right)S_{i,j}-\frac{1}{2w}\left(D_iF_j+D_jF_i\right)+\frac{1}{2w^2}D_iD_j\right\}x_ix_j\\
&+\sum_{i=1}^{n}\left(\frac{\partial N_i}{\partial s}-\lambda_iN_i+M_1F_i-\frac{1}{w}\left(\sum_{k=1}^{n}c_kB_k\right)F_i-\frac{1}{w}\left(\sum_{k=1}^{n}c_kN_k\right)D_i+\frac{1}{w^2}\left(\sum_{k=1}^{n}c_kB_k\right)D_i\right)x_i\,,\\
&+\frac{\partial\tilde{E}}{\partial s}-C+\frac{1}{2}M_2\sum_{i=1}^{n}c_iS_{i,i}+M_1\sum_{i=1}^{n}c_iN_i-\frac{1}{w}\left(\sum_{k=1}^{n}c_kB_k\right)\sum_{i=1}^{n}c_iN_i+\frac{1}{2w^2}\left(\sum_{i=1}^{n}c_iB_i\right)=0
\end{aligned}
\tag{84}
$$

leading to a system

$$
0=\frac{\partial S_{i,j}}{\partial s}-\left(\lambda_i+\lambda_j\right)S_{i,j}-\frac{1}{w}\left(D_iF_j+D_jF_i\right)+\frac{1}{w^2}D_iD_j\quad(1\le i,j\le n,\ \ 0\le s\le P),\tag{85}
$$

$$
\begin{aligned}
0=&\frac{\partial N_i}{\partial s}-\lambda_iN_i+M_1F_i\\
&-\frac{1}{w}\left(\sum_{k=1}^{n}c_kB_k\right)F_i-\frac{1}{w}\left(\sum_{k=1}^{n}c_kN_k\right)D_i+\frac{1}{w^2}\left(\sum_{k=1}^{n}c_kB_k\right)D_i
\end{aligned}
\quad(1\le i\le n,\ \ 0\le s\le P),\tag{86}
$$

$$
0=\frac{\partial\tilde{E}}{\partial s}-C+M_1\sum_{i=1}^{n}c_iN_i+\frac{1}{2}M_2\sum_{i=1}^{n}c_iS_{i,i}-\frac{1}{w}\left(\sum_{i=1}^{n}c_iB_i\right)\left(\sum_{i=1}^{n}c_iN_i\right)+\frac{1}{2w^2}\left(\sum_{i=1}^{n}c_iB_i\right)\quad(0\le s\le P).\tag{87}
$$

Finally, we obtain $C$ as

$$
C=\frac{1}{T}\int_0^T\left\{\frac{1}{2}M_2\sum_{i=1}^{n}c_iS_{i,i}+M_1\sum_{i=1}^{n}c_iN_i-\frac{1}{w}\left(\sum_{k=1}^{n}c_kB_k\right)\sum_{i=1}^{n}c_iN_i+\frac{1}{2w^2}\left(\sum_{i=1}^{n}c_iB_i\right)\right\}\mathrm{d}s\,.\tag{88}
$$

### 3.5 Performance evaluation in an infinite-dimensional case

We also present an infinite-dimensional counterpart of the KBE: for $Y\in L^2(\pi)$ and $0\le s\le P$,

$$
\begin{aligned}
&-C+\frac{\partial W}{\partial s}+\int_0^{+\infty}\nabla W(Y)\left(-\lambda Y(\lambda)-\frac{1}{w}\left(\int_0^{+\infty}\int_0^{+\infty}\Gamma(\omega,\theta)Y(\theta)\pi(\mathrm{d}\theta)+\gamma(\omega)\right)\pi(\mathrm{d}\omega)\right)\pi(\mathrm{d}\lambda)\\
&+\int_0^{+\infty}\left(\int_0^{+\infty}W\left(Y(\cdot)+z\hat{\delta}_{\{\omega=(\cdot)\}}\right)\pi(\mathrm{d}\omega)-W(Y)\right)v(\mathrm{dz})\\
&+\frac{1}{2}\left(-\frac{1}{w}\int_0^{+\infty}\left(\int_0^{+\infty}\Gamma(\lambda,\theta)Y(\theta)\pi(\mathrm{d}\theta)+\gamma(\lambda)\right)\pi(\mathrm{d}\lambda)\right)^2=0
\end{aligned}
\tag{89}
$$

whose solution is the couple $(W,C)$ of $W:[0,P]\times L^2(\pi)\to\mathbb{R}$ with the period $P$ and $C\in\mathbb{R}$, where the dependence of variables on time is suppressed. Using an ansatz with some symmetric and time-dependent $\Lambda(\cdot,\cdot)$, a time-dependent field $\rho(\cdot)$, and a time-dependent scalar $G$:

$$
W(Y)=\frac{1}{2}\int_0^{+\infty}\int_0^{+\infty}\Lambda(\theta,\lambda)Y(\lambda)Y(\theta)\pi(\mathrm{d}\lambda)\pi(\mathrm{d}\theta)+\int_0^{+\infty}\rho(\lambda)Y(\lambda)\pi(\mathrm{d}\lambda)+G\,.\tag{90}
$$

After a straightforward calculation using (89)–(90), we obtain

$$
\begin{aligned}
0 &= -C + \frac{\partial W}{\partial s} + M_1 \int_0^{+\infty} \rho(\lambda)\pi(\mathrm{d}\lambda) + \frac{1}{2} M_2 \int_0^{+\infty} \Lambda(\theta,\theta)\pi(\mathrm{d}\theta) \\
&- \frac{1}{w}\left(\int_0^{+\infty} \gamma(\lambda)\pi(\mathrm{d}\lambda)\right)\left(\int_0^{+\infty} \rho(\lambda)\pi(\mathrm{d}\lambda)\right) + \frac{1}{2w^2}\left(\int_0^{+\infty} \gamma(\lambda)\pi(\mathrm{d}\lambda)\right)^2 \\
&+ \int_0^{+\infty} \left(\begin{array}{l} -\lambda\rho(\lambda) + M_1 \int_0^{+\infty} \Lambda(\theta,\lambda)\pi(\mathrm{d}\theta) \\ -\frac{1}{w}\int_0^{+\infty}\left(\left(\int_0^{+\infty}\gamma(\omega)\pi(\mathrm{d}\omega)\right)\Lambda(\lambda,\theta) + \left(\int_0^{+\infty}\rho(\omega)\pi(\mathrm{d}\omega)\right)\Gamma(\lambda,\theta)\right)\pi(\mathrm{d}\theta) \\ +\frac{1}{w^2}\int_0^{+\infty}\gamma(\theta)\pi(\mathrm{d}\theta)\int_0^{+\infty}\Gamma(\lambda,\theta)\pi(\mathrm{d}\theta) \end{array}\right) Y(\lambda)\pi(\mathrm{d}\lambda) \\
&+ \int_0^{+\infty}\int_0^{+\infty} \left\{\begin{array}{l} -\frac{1}{2}(\theta+\lambda)\Lambda(\theta,\lambda) \\ -\frac{1}{2w}\int_0^{+\infty}\int_0^{+\infty}\left(\Gamma(\lambda,\omega)\Lambda(\theta,\sigma) + \Gamma(\theta,\omega)\Lambda(\lambda,\sigma)\right)\pi(\mathrm{d}\sigma)\pi(\mathrm{d}\omega) \\ +\frac{1}{2w^2}\left(\int_0^{+\infty}\Gamma(\lambda,\theta_1)\pi(\mathrm{d}\theta_1)\right)\left(\int_0^{+\infty}\Gamma(\lambda_1,\theta)\pi(\mathrm{d}\lambda_1)\right) \end{array}\right\} Y(\theta)Y(\lambda)\lambda\pi(\mathrm{d}\theta)\pi(\mathrm{d}\lambda)
\end{aligned}
\tag{91}
$$

We then get partial integro-differential equations for $\lambda,\theta > 0$, $s \in [0,P]$ associated with (90):

$$
\begin{aligned}
0 &= \frac{\partial \Lambda(\theta,\lambda)}{\partial s} - (\theta+\lambda)\Lambda(\theta,\lambda) \\
&- \frac{1}{w}\int_0^{+\infty}\int_0^{+\infty}\left(\Gamma(\lambda,\omega)\Lambda(\theta,\sigma) + \Gamma(\theta,\omega)\Lambda(\lambda,\sigma)\right)\pi(\mathrm{d}\sigma)\pi(\mathrm{d}\omega), \\
&+ \frac{1}{w^2}\left(\int_0^{+\infty}\Gamma(\lambda,\omega)\pi(\mathrm{d}\omega)\right)\left(\int_0^{+\infty}\Gamma(\omega,\theta)\pi(\mathrm{d}\omega)\right)
\end{aligned}
\tag{92}
$$

$$
\begin{aligned}
0 &= \frac{\partial \rho(\lambda)}{\partial s} - \lambda\rho(\lambda) + M_1\int_0^{+\infty}\Lambda(\theta,\lambda)\pi(\mathrm{d}\theta) \\
&- \frac{1}{w}\int_0^{+\infty}\left(\left(\int_0^{+\infty}\gamma(\omega)\pi(\mathrm{d}\omega)\right)\Lambda(\lambda,\theta) + \left(\int_0^{+\infty}\rho(\omega)\pi(\mathrm{d}\omega)\right)\Gamma(\lambda,\theta)\right)\pi(\mathrm{d}\theta), \\
&+ \frac{1}{w^2}\int_0^{+\infty}\gamma(\theta)\pi(\mathrm{d}\theta)\int_0^{+\infty}\Gamma(\lambda,\theta)\pi(\mathrm{d}\theta)
\end{aligned}
\tag{93}
$$

$$
\begin{aligned}
0 &= \frac{\partial G}{\partial s} - C + M_1\int_0^{+\infty}\rho(\lambda)\pi(\mathrm{d}\lambda) + \frac{1}{2}M_2\int_0^{+\infty}\Lambda(\theta,\theta)\pi(\mathrm{d}\theta) \\
&- \frac{1}{w}\left(\int_0^{+\infty}\gamma(\lambda)\pi(\mathrm{d}\lambda)\right)\left(\int_0^{+\infty}\rho(\lambda)\pi(\mathrm{d}\lambda)\right) + \frac{1}{2w^2}\left(\int_0^{+\infty}\gamma(\lambda)\pi(\mathrm{d}\lambda)\right)^2.
\end{aligned}
\tag{94}
$$

**Remark 6**: The discretized counterparts of (92)–(94) are (85)–(88). As in **Remark 5**, they are consistent with each other.

### 3.6 Optimality

The optimality of the HJB and Riccati equations is verified. The LQ formulation allows for analysis in the context of smooth solutions. Our strategy is based on Itô's formula for infinite-dimensional stochastic processes (e.g., Kadlec and Maslowski, 2019). The main result here is the following **Proposition 3**, whose proof is presented in **Appendix C.3**.

***Proposition 3***

*Assume that a smooth couple* $\left(\Gamma_t(\cdot,\cdot),\gamma_t(\cdot)\right)_{0\le t\le P}$ *with a positive semi-definite* $\Gamma_t : L^{\infty}(\pi\otimes\pi)\to\mathbb{R}$ *and* $\gamma_t : L^{\infty}(\pi)\to\mathbb{R}$, *both having the period* $P$ *and continuously differentiable for* $0\le t\le P$, *satisfies (67)-(68). If the control* $u^*$ *of (71) is admissible and the controlled process* $(Y_t)_{0\le t\le P}$ *is generated by a periodic probability measure, then* $J$ *is minimized by* $u=u^*$. *We then have* $J(u^*)=H$ *with* $H$ *in (69).*

We also state a linkage between the HJB and Riccati equations.

***Definition 1***

*A classical solution to the HJB equation (57) is a couple* $(V,H)$ *of* $V:[0,P]\times L^2(\pi)\to\mathbb{R}$ *and* $H\in\mathbb{R}$. *Here, at each* $t\in[0,P]$, $V$ *is continuously Fréchet differentiable in* $L^2(\pi)$ *with* $\nabla V(Y)$ *belongings to* $L^2(\pi)$ *given* $Y\in L^2(\pi)$, $V$ *is continuously differentiable with respect to* $t\in[0,P]$, *and* $(V,H)$ *satisfies (57).*

We show that the couple $(V,H)$ of $V$ given by (58) and $H$ obtained from (69) is a solution in the sense of **Definition 1**. The following **Proposition 4** states that the HJB equation derived heuristically is the optimality equation of the LQ control problem. This is a direct consequence of **Proposition 3**.

***Proposition 4***

*Under the assumption of* ***Proposition 3***, $V$ *given by (58) and* $H$ *of (69) is a solution* $(V,H)$ *to the HJB equation (57).*

## 4 Application

### 4.1 Target site

The target sites in this work include the three observation stations (Stations Y, D, U) of the Hiikawa River in Japan around Obara Dam (called Station D) controlling midstream to downstream reaches of the river (**Figure 1**). A public record of the inflow from an upstream river of the dam is also available (http://www1.river.go.jp/cgi-bin/SrchDamData.exe?ID=607041287705020&KIND=1&PAGE=0) provided by Ministry of Land, Infrastructure, Transport and Tourism (MLIT)), which is called Station U for brevity. We measured the water depth and temperature (10-min interval: HOBO U20) at Station Y 8.2 km downstream of the dam. Station Y is chosen because there is a good fishing ground for the amphidromous fish *Plecoglossus altivelis*, serving as one of the most important commercial fishery resources in Japan. Its growth has been considered to be negatively affected by the growth of benthic algae due to regulated streamflow (Yoshioka and Yaegashi, 2018). Station Y is located near a mountainous hot spring area called Yumura Onsen, which is often visited by tourists. Preserving the river environment around this station is therefore important from both ecological and economic perspectives. Station U is considered because the probability density of discharge at this station is qualitatively different from those at Stations Y and D. Namely, we consider Station U to demonstrate that the supOU process can be used for a variety of discharge data in the study area.

### 4.2 Data curation

The public discharge data of MLIT is used for Station D and U. Discharge data at Station Y is not available and is therefore emulated by a longitudinally one-dimensional hydraulic model (Tanaka et al., 2021) as explained in **Appendix A**. This model is based on a local inertial equation as a widely used model for simulating flood events (de Almeida and Bates, 2013). For all the stations, the data for analysis is from June 1, 2016, to May 31, 2020. **Figure 2** presents discharge times series at each station. We set $P$ as one year (365.25 days).

### 4.3 Model identification

Hourly discharge data was used for the identification of the supOU process at each station. The identification method of Yoshioka and Yoshioka (2021) was extended so that all parameters were identified. The identification procedure consisted of two steps. The first step was to identify the parameters $B_\pi$ and $\alpha_\pi$ of $\pi$ by a least-squares fitting between the theoretical and empirical ACFs. After identifying $B_\pi$ and $\alpha_\pi$, the second step was to identify the minimum discharge $\underline{X}$ and parameters $a_v$, $b_v$, $p_v$, $\alpha_v$ of $v$. The minimum discharge $\underline{X}$ was found from the time series. The remaining parameters were identified as the minimizer of the following metric:

$$\left(\frac{\mathrm{Ave}_{\mathrm{model}}-\mathrm{Ave}_{\mathrm{Data}}}{\mathrm{Ave}_{\mathrm{Data}}}\right)^2+\left(\frac{\mathrm{Std}_{\mathrm{model}}-\mathrm{Std}_{\mathrm{Data}}}{\mathrm{Std}_{\mathrm{Data}}}\right)^2+\left(\frac{\mathrm{Skew}_{\mathrm{model}}-\mathrm{Skew}_{\mathrm{Data}}}{\mathrm{Skew}_{\mathrm{Data}}}\right)^2+\left(\frac{\mathrm{Kurt}_{\mathrm{model}}-\mathrm{Kurt}_{\mathrm{Data}}}{\mathrm{Kurt}_{\mathrm{Data}}}\right)^2, \tag{95}$$

where the subscripts "Model" and "Data" classify the modelled ones based on the supOU process (7)–(9) and the data at each station.

The identified model was examined against the data from the viewpoint of probability density. The probability density of discharge is not available in a closed form, which was therefore computed by a Monte-Carlo (MC) method (Yoshioka, 2022). This numerical method is based on the Euler–Maruyama scheme and observation that the reversion speed is generated randomly from $\pi$ at each jump (Fasen and Klüppelberg, 2007; Stelzer et al., 2015):

$$X_{(n+1)\Delta t}=e^{-\lambda_n\Delta t}X_{n\Delta t}+\left(1-e^{-\lambda_n\Delta t}\right)\underline{X}+\frac{1-e^{-\lambda_n\Delta t}}{\lambda_n\Delta t}\left(\Delta L_n\right)^{\frac{1}{p_v}},\quad n=0,1,2,\ldots \tag{96}$$

with $X_0>0$. Here, $\Delta t\left(=0.001\right)$ (h) is the time increment, each $\lambda_n$ is a random variable generated by the identified density $\pi$, each $\Delta L_n$ is a tempered stable process with the time increment $\Delta t$ generated by the rejection sampling method (Kawai and Masuda, 2011), and each $\Delta L_n$ is mutually independent and is associated with the Lévy measure $v'$ given by

$$v'(\mathrm{dz}) = \frac{a_v}{p_v} \frac{e^{-(b_v)^{p_v} z}}{z^{1+\frac{\alpha_v}{p_v}}} \mathrm{dz}, \quad z > 0. \tag{97}$$

Here, we invoked an equivalence between the tempered stable and generalized tempered stable subordinators; for each $k \geq 1$, using (6), we obtain

$$\int_0^{+\infty} (z)^{\frac{k}{p_v}} v'(\mathrm{dz}) = \int_0^{+\infty} \frac{a_v}{p_v} (z)^{\frac{k}{p_v} - \frac{\alpha_v}{p_v} - 1} e^{-(b_v)^{p_v} z} \mathrm{dz} = \frac{a_v}{p_v} (b_v)^{\frac{\alpha_v - k}{p_v}} \Gamma\left(\frac{k - \alpha_v}{p_v}\right) = M_k. \tag{98}$$

This MC method applies only to the uncontrolled case as it cannot handle state-dependent drift.

**Tables 1** and **2** present the identified model parameters and the corresponding statistics of the identified model and empirical ones with their relative errors, respectively. **Figure 3** shows the empirical and fitted ACFs at each observation station, exhibiting that the proposed model accurately captures the sub-exponential decay. **Table 2** demonstrates that, for each observation station, the identified model with the parameters of **Table 1** successfully reproduces the statistical moments. The supporting data presented in **Appendix B** implies that the generalized tempered stable ($p_v = 2$) model more accurately reproduces the statistics than the tempered stable model ($p_v = 1$). **Figure 4** compares the empirical and modeled probability density functions (PDFs) at each observation station where the latter is generated by (96). The identified model captures the unimodal nature of the PDFs at Stations Y and D and the decreasing nature of the PDF at Station U. Furthermore, **Table 3** compares the theoretical and computed statistics, demonstrating their good agreement and thus the validity of the MC simulation.

In summary, the identification and verification results suggest that the supOU process serves as a viable mathematical tool for analyzing discharge time series in perennial river environments.

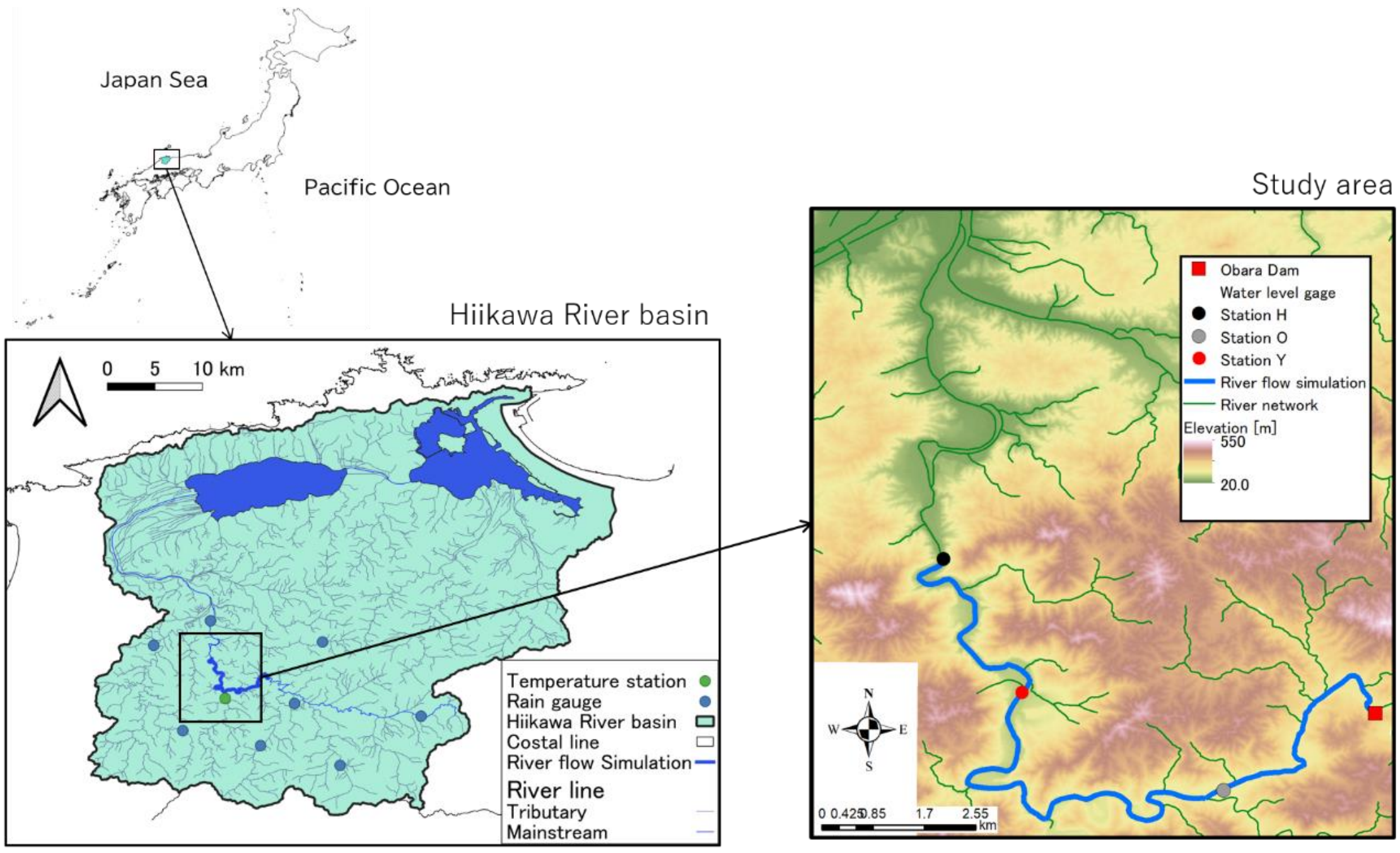

**Figure 1.** Map of the study area. Stations H and O are not referred to in the main text but used in the development of the hydraulic model in **Appendix A**.

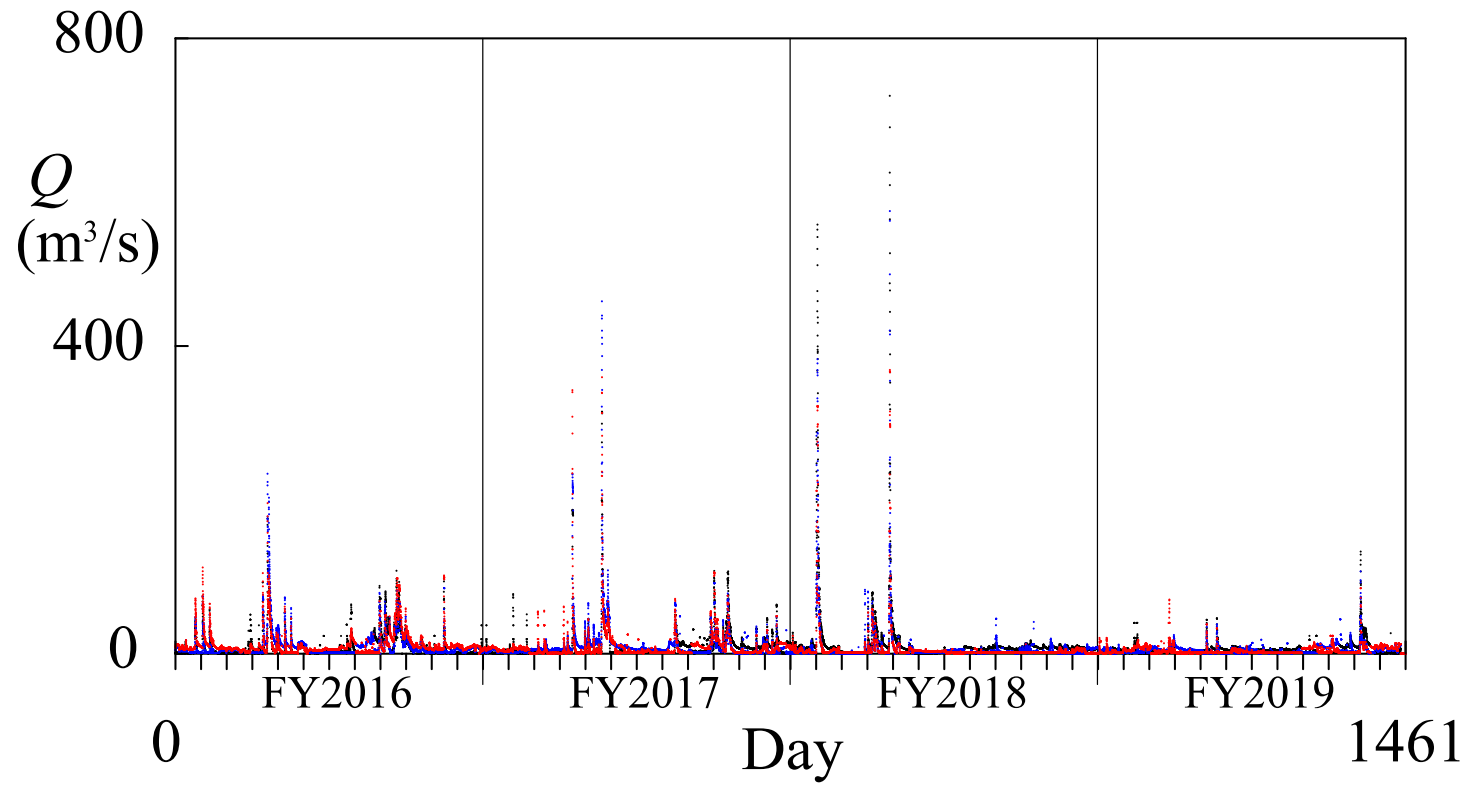

**Figure 2.** Discharge data at each station, Station Y (Red), Station D (Blue), Station U (Black).

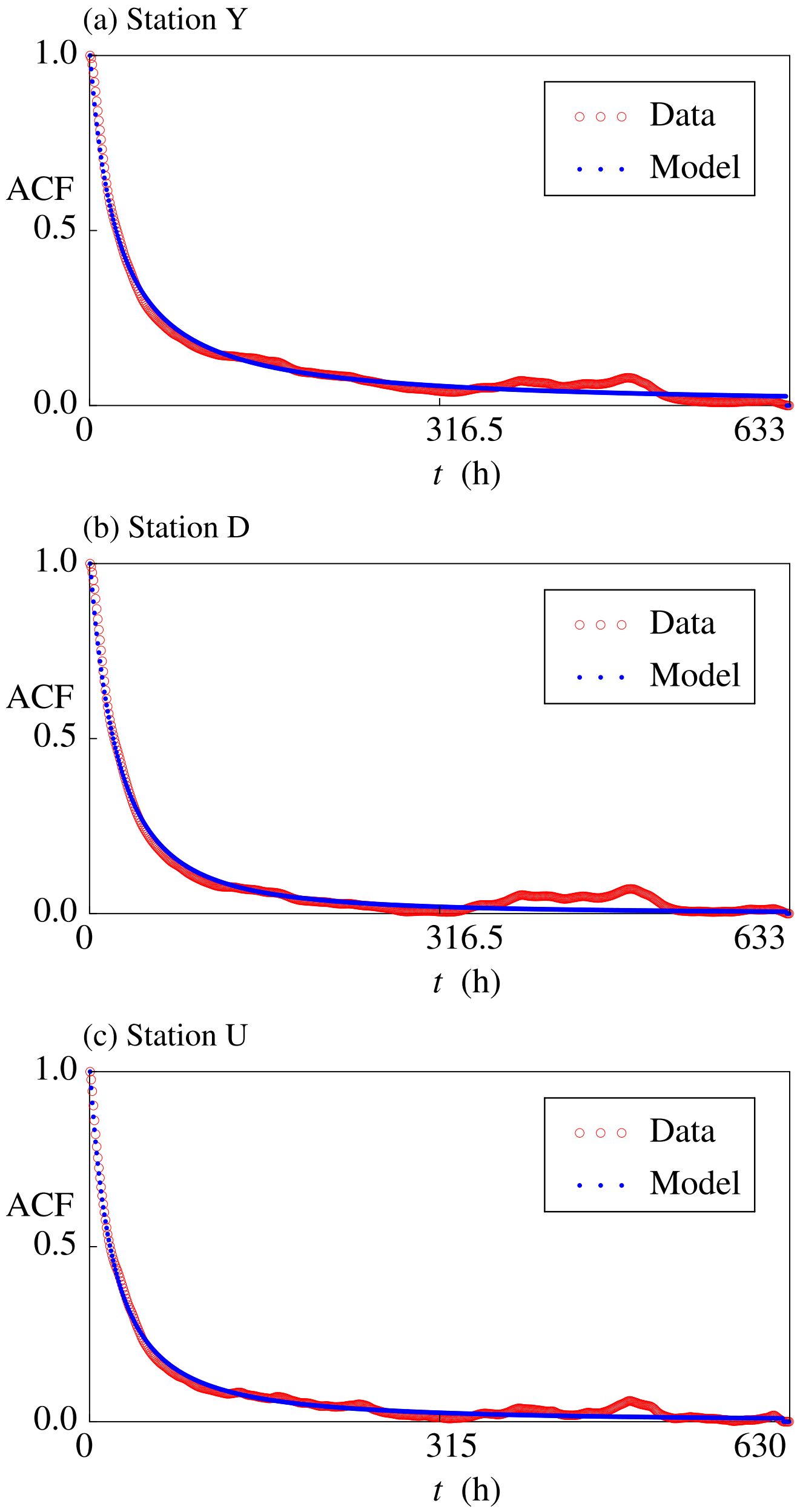

**Figure 3.** Theoretical and empirical ACFs at each station: (a) Station Y, (b) Station D, (c) Station U.

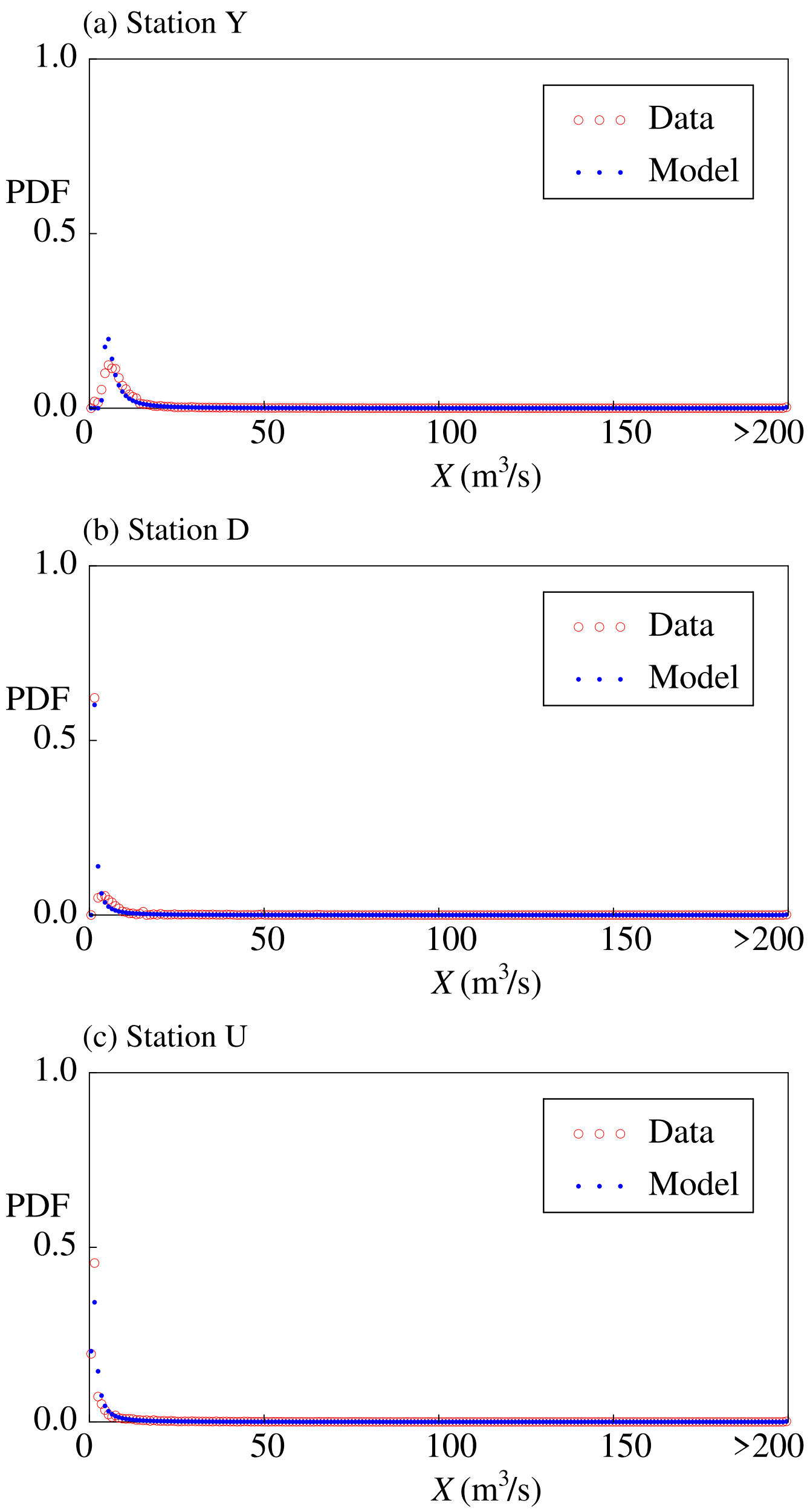

**Figure 4.** Theoretical and empirical PDFs at each station (a) Station Y, (b) Station D, (c) Station U.

**Table 1.** Identified parameter values.

| Parameter | Station Y | Station D | Station U |
|---|---|---|---|
| $B_\pi$ (1/h) | 0.0344 | 0.0201 | 0.0315 |
| $\alpha_\pi$ (-) | 2.17 | 2.97 | 2.53 |
| $\underline{X}$ ($m^3/s$) | 1.28 | 1.00 | 0.00 |
| $p_v$ (-) | 2 | 2 | 2 |
| $\alpha_v$ (-) | 4.08.E-01 | 5.25.E-01 | 7.05.E-01 |
| $a_v$ ($m^{3\alpha_v}/s^{\alpha_v}$) | 1.27.E-02 | 5.18.E-03 | 1.34.E-02 |
| $b_v$ ($s/m^3$) | 2.33.E-06 | 7.73.E-06 | 4.59.E-06 |

**Table 2.** Statistics of the identified model and empirical ones.

| | Station Y | | Station D | | Station U | |
|---|---|---|---|---|---|---|
| Least-squares error | 0.00218 | | 0.00238 | | 0.000107 | |
| | Data | Model | Data | Model | Data | Model |
| Ave ($m^3/s$) | 1.21.E+01 | 1.21.E+01 | 5.13.E+00 | 5.11.E+00 | 5.40.E+00 | 5.39.E+00 |
| Std ($m^3/s$) | 2.26.E+01 | 2.27.E+01 | 1.54.E+01 | 1.55.E+01 | 1.66.E+01 | 1.66.E+01 |
| Skew (-) | 1.28.E+01 | 1.22.E+01 | 1.19.E+01 | 1.13.E+01 | 1.27.E+01 | 1.26.E+01 |
| Kurt (-) | 2.43.E+02 | 2.48.E+02 | 1.95.E+02 | 1.99.E+02 | 2.54.E+02 | 2.55.E+02 |

**Table 3.** Comparison of statistics between the theory and MC simulation.

| | Station Y | | Station D | | Station U | |
|---|---|---|---|---|---|---|
| | Theory | MC | Theory | MC | Theory | MC |
| Ave ($m^3/s$) | 1.21.E+01 | 1.17.E+01 | 5.11.E+00 | 5.11.E+00 | 5.39.E+00 | 5.39.E+00 |
| Std ($m^3/s$) | 2.27.E+01 | 2.26.E+01 | 1.55.E+01 | 1.55.E+01 | 1.66.E+01 | 1.67.E+01 |
| Skew (-) | 1.22.E+01 | 1.19.E+01 | 1.13.E+01 | 1.12.E+01 | 1.26.E+01 | 1.25.E+01 |
| Kurt (-) | 2.48.E+02 | 2.46.E+02 | 1.99.E+02 | 1.97.E+02 | 2.55.E+02 | 2.56.E+02 |

### 4.4 Computational method

The finite-dimensional Riccati equation is discretized in time by a classical forward Euler method with time increment $\Delta t = O\left(n^{-1}\right) > 0$. A higher-order method such as a Runge–Kutta method can also be used, but the Euler method will turn out to be sufficient for our case due to the dominance of the spatial discretization error for relatively fine computational resolution. We must take $\Delta t$ to be sufficiently small because the forward Euler method is explicit in time.

We cyclically solve the Riccati equation backward in time starting from a zero initial guess for several periods until the numerical solution converges sufficiently. The convergence speed is expected to be exponential according to Theorem 9 of Ha and Choi (2019). In our case, integrating the equation for less than 10 periods was sufficient for the convergence with the maximum absolute difference smaller than $10^{-10}$ in the latest successive two periods. We set $\Delta t = 1/(24n)$ (day). The following discretization parameters were used unless otherwise specified: $n = 160$, $\bar{\eta} = 0.02$ (1/day).

### 4.5 Verification with a manufactured solution

We verify the convergence of the numerical method against a manufactured solution to an exactly solvable system. We consider a periodic solution to the following system for $\lambda, \theta > 0$, $s \in [0, P]$:

$$0 = \frac{\partial \Gamma_s(\theta,\lambda)}{\partial s} - (\theta+\lambda)\Gamma_s(\theta,\lambda) - \frac{1}{w}\int_0^{+\infty}\Gamma_s(\theta_1,\lambda)\pi(\mathrm{d}\theta_1)\int_0^{+\infty}\Gamma_s(\theta_2,\theta)\pi(\mathrm{d}\theta_2) + 1 + f_\Gamma(s,\theta,\lambda), \quad (99)$$

$$\begin{aligned} 0 &= \frac{\partial \gamma_s(\lambda)}{\partial s} - \lambda\gamma_s(\lambda) \\ &\quad - \frac{1}{w}\int_0^{+\infty}\Gamma_s(\theta,\lambda)\pi(\mathrm{d}\theta)\int_0^{+\infty}\gamma_s(\lambda)\pi(\mathrm{d}\lambda) + M_1\int_0^{+\infty}\Gamma_s(\theta,\lambda)\pi(\mathrm{d}\theta) - \bar{X}_s + g_\Gamma(s,\lambda) \end{aligned}, \quad (100)$$

$$\begin{aligned} 0 &= \frac{\partial E}{\partial s} - H - \frac{1}{2w}\left(\int_0^{+\infty}\gamma_s(\lambda)\pi(\mathrm{d}\lambda)\right)^2 \\ &\quad + \frac{1}{2}M_2\int_0^{+\infty}\Gamma_s(\theta,\theta)\pi(\mathrm{d}\theta) + M_1\int_0^{+\infty}\gamma_s(\lambda)\pi(\mathrm{d}\lambda) + \frac{1}{2}\left(\bar{X}_s\right)^2 \end{aligned}. \quad (101)$$

The artificial sources $f_\Gamma$ and $g_\Gamma$ are chosen so that (99)–(100) admit the smooth solution

$$\Gamma_s(\lambda,\theta) = a_\Gamma(s)\exp\left(-b_\Gamma(\lambda+\theta)\right) \text{ and } \gamma_s(\lambda) = c_\Gamma(s)\exp\left(-b_\Gamma\lambda\right) \quad (102)$$

with smooth coefficients $a_\Gamma, c_\Gamma > 0$ and a constant $b_\Gamma > 0$. We should set $f_\Gamma$ and $g_\Gamma$ as

$$f_\Gamma(s,\theta,\lambda) = -1 + \left(-\frac{\mathrm{d}a_\Gamma}{\mathrm{d}s} + (\theta+\lambda)a_\Gamma + \frac{(a_\Gamma)^2}{w\left(B_\pi b_\Gamma + 1\right)^{2\alpha_\pi}}\right)\exp\left(-b_\Gamma(\lambda+\theta)\right), \quad (103)$$

$$g_\Gamma(s,\lambda) = \bar{X}_s + \left(-\frac{\mathrm{d}c_\Gamma}{\mathrm{d}s} + \lambda c_\Gamma + \frac{a_\Gamma c_\Gamma}{w\left(B_\pi b_\Gamma + 1\right)^{2\alpha_\pi}} - \frac{M_1 a_\Gamma}{\left(B_\pi b_\Gamma + 1\right)^{\alpha_\pi}}\right)\exp\left(-b_\Gamma\lambda\right). \quad (104)$$

We then obtain $H$ as

$$H = \frac{1}{P}\int_0^P \left( -\frac{\left(c_\Gamma(s)\right)^2}{2w\left(B_\pi b_\Gamma + 1\right)^{2\alpha_\pi}} + \frac{M_2 a_\Gamma(s)}{2\left(2B_\pi b_\Gamma + 1\right)^{\alpha_\pi}} + \frac{M_1 c_\Gamma(s)}{\left(B_\pi b_\Gamma + 1\right)^{\alpha_\pi}} + \frac{\left(\bar{X}_s\right)^2}{2} \right) \mathrm{d}s \,. \tag{105}$$

The modified equations share the same integral operators with the original Riccati equation, suggesting that the convergence of the numerical method can be examined with the manufactured solution.

We set $a_\Gamma = 0.1 + 0.05\sin\left(\frac{2\pi s}{P}\right)$, $c_\Gamma = 0.2 + 0.1\sin\left(\frac{2\pi s}{P}\right)$, and $b_\Gamma = 0.02$ (h). In both cases, $H$ is evaluated analytically. We use the identified model at Station Y. The target discharge is set as $\hat{X}_s = 10\left(1 + 0.5\cos\left(\frac{2\pi s}{P}\right)\right)$ and $w = 1$. We examine different values of $(n, \beta)$ to verify the convergence of the proposed numerical method.

**Tables 4–6** present the manufactured and computed $H$ and the corresponding relative errors given by $e_n = \frac{\left|H_{\mathrm{Manufactured}} - H_{\mathrm{Computed}}\right|}{H_{\mathrm{Manufactured}}}$ for $\beta = 0.2,\ 0.5,\ 0.8$ , respectively. Here, $H_{\mathrm{Manufactured}} = 46.2495$ ($\mathrm{m}^6/\mathrm{s}^2$). The convergence rates between each discretization level, given by $c_n = \log_{\frac{n_1}{n}}\left(\frac{e_n}{e_{n_1}}\right)$ ($n_1 > n$, $n$: current discretization level, $n_1$: next discretization level), are also provided in the tables. The computational results suggest that the relative error decreases as the discretization level increases.

We also compare the manufactured and computed $\Gamma$, $\gamma$ to analyze the convergence of numerical solutions more in detail. **Tables 7–9** present the $l^\infty$ errors, maximal errors among at all space-time points $\left(\eta_i^n, \eta_j^n, k\Delta t\right)$ in $\left(0, \eta_{n,n}\right)^2 \times \left[0, P\right]$, and corresponding convergence rates between the manufactured and computed $\Gamma$, $\gamma$. The convergence order is $O(1)$ for $\beta = 0.2,\ 0.5$, which is not prohibitive as standard modern numerical solvers exhibit a convergence order around 1. It also suggests that the use of the forward Euler method is sufficient in these cases because the order of convergence is close to 1 with respect to $n$. By contrast, numerical solutions do not converge for $\beta = 0.8$. The computational results suggest that choosing too large $\beta$ with which the domain truncation is relatively strong should be avoided for efficiently computing the Riccati equation.

Although our results are based on computational experiments, they imply that the proposed numerical method is convergent if the discretization parameters are specified properly.

**Table 4.** Relative error and convergence rate for $\beta = 0.2$ .

| Discretization level $n$ | $H_{\text{Computed}}$ | $e_n$ | $r_{n,n_1}$ |
|---|---|---|---|
| 10 | 46.1774 | 7.21.E-02 | 3.54.E+00 |
| 20 | 46.2433 | 1.34.E-04 | 6.55.E+00 |
| 40 | 46.2495 | 1.43.E-06 | 3.91.E+00 |
| 80 | 46.2495 | 9.52.E-08 | 4.17.E-01 |
| 160 | 46.2495 | 7.13.E-08 | 4.20.E-01 |
| 320 | 46.2495 | 5.33.E-08 | |

**Table 5.** Relative error and convergence rate for $\beta = 0.5$.

| Discretization level $n$ | $H_{\text{Computed}}$ | $e_n$ | $r_{n,n_1}$ |
|---|---|---|---|
| 10 | 45.9120 | 7.30.E-03 | 9.57.E-01 |
| 20 | 46.0756 | 3.76.E-03 | 1.35.E+00 |
| 40 | 46.1811 | 1.48.E-03 | 1.96.E+00 |
| 80 | 46.2319 | 3.80.E-04 | 2.90.E+00 |
| 160 | 46.2472 | 5.08.E-05 | 4.33.E+00 |
| 320 | 46.2494 | 2.52.E-06 | |

**Table 6.** Relative error and convergence rate for $\beta = 0.8$.

| Discretization level $n$ | $H_{\text{Computed}}$ | $e_n$ | $r_{n,n_1}$ |
|---|---|---|---|
| 10 | 45.4822 | 1.66.E-02 | 2.01.E-01 |
| 20 | 45.5821 | 1.44.E-02 | 2.18.E-01 |
| 40 | 45.6758 | 1.24.E-02 | 2.41.E-01 |
| 80 | 45.7640 | 1.05.E-02 | 2.68.E-01 |
| 160 | 45.8464 | 8.72.E-03 | 3.02.E-01 |
| 320 | 45.9225 | 7.07.E-03 | |

**Table 7.** $l^\infty$ error of computed $\Gamma$, $\gamma$ for $\beta = 0.2$ .

| | Error | | Convergence rate | |
|---|---|---|---|---|
| Discretization level $n$ | $\Gamma$ | $\gamma$ | $\Gamma$ | $\gamma$ |
| 10 | 3.90.E-02 | 1.86.E-01 | 3.20.E+00 | 3.06.E+00 |
| 20 | 4.25.E-03 | 2.24.E-02 | 5.95.E+00 | 5.91.E+00 |
| 40 | 6.86.E-05 | 3.73.E-04 | 3.90.E+00 | 5.27.E+00 |
| 80 | 4.60.E-06 | 9.64.E-06 | 9.89.E-01 | 9.70.E-01 |
| 160 | 2.32.E-06 | 4.92.E-06 | 9.82.E-01 | 9.48.E-01 |
| 320 | 1.17.E-06 | 2.55.E-06 | | |

**Table 8.** $l^\infty$ error of computed $\Gamma$, $\gamma$ for $\beta = 0.5$.

| | Error | | Convergence rate | |
|---|---|---|---|---|
| Discretization level $n$ | $\Gamma$ | $\gamma$ | $\Gamma$ | $\gamma$ |
| 10 | 2.86.E-01 | 1.57.E+00 | 6.84.E-01 | 3.39.E-01 |
| 20 | 1.78.E-01 | 1.24.E+00 | 9.38.E-01 | 7.35.E-01 |
| 40 | 9.30.E-02 | 7.44.E-01 | 1.43.E+00 | 1.32.E+00 |
| 80 | 3.44.E-02 | 2.98.E-01 | 2.27.E+00 | 2.20.E+00 |
| 160 | 7.17.E-03 | 6.49.E-02 | 3.53.E+00 | 3.49.E+00 |
| 320 | 6.21.E-04 | 5.76.E-03 | | |

**Table 9.** $l^\infty$ error of computed $\Gamma$, $\gamma$ for $\beta = 0.8$.

| | Error | | Convergence rate | |
|---|---|---|---|---|
| Discretization level *n* | $\Gamma$ | $\gamma$ | $\Gamma$ | $\gamma$ |
| 10 | 1.35.E+00 | 6.29.E+00 | -2.96.E-01 | -6.74.E-01 |
| 20 | 1.66.E+00 | 1.00.E+01 | -3.25.E-01 | -5.52.E-01 |
| 40 | 2.08.E+00 | 1.47.E+01 | -3.43.E-01 | -4.69.E-01 |
| 80 | 2.63.E+00 | 2.04.E+01 | -3.58.E-01 | -4.53.E-01 |
| 160 | 3.38.E+00 | 2.79.E+01 | -4.13.E-01 | -5.01.E-01 |
| 320 | 4.50.E+00 | 3.95.E+01 | | |

### 4.5 Application to algae bloom control

#### 4.5.1 Computational setting

The computation here focuses on an effective suppression of the thick growth of benthic algae in river environments (Huang et al., 2021; Loire et al., 2021; Yoshioka and Tsujimura, 2020) by installing some infrastructure that can add/subtract water at the target site. Benthic algae such as *Cladophora glomerata (C. glomerata)* grow during the summer high-temperature season especially if the discharge, namely hydraulic force to remove them from the riverbed, is not sufficiently large. In our context, the discharge should be controlled so that it is sufficiently large during summer with the least cost. The alga considered here is *C. glomerata,* whose growth dynamics have been studied extensively. Water temperature is a key factor that affects the growth of *C. glomerata*. It has been reported to persist between the water temperature above around 5 °C and below around 24 °C with an optimum of around 13–17 °C (Graham et al., 1982; Dodds, 1991). Yoshioka and Yaegashi (2018) found a critical water discharge around 20 (m$^3$/s) in this river above which the growth of *C. glomerata* is limited.

Based on these observations and **Remark 1**, we set $J$ as

$$J = \limsup_{k \to +\infty} \frac{1}{kP} \mathbb{E}\left[ \int_0^{kP} \frac{w'_s}{2} \left( X_{n,s} - \hat{X}_s \right)^2 + w \frac{u_s^2}{2} \mathrm{d}s \right] \tag{106}$$

with $\hat{X}_s = 20$ (m$^3$/s) and $w'_s = \varepsilon + 4\left(\bar{W} - \underline{W}\right)^{-2} \max\left\{ \left(\bar{W} - W_s\right)\left(W_s - \underline{W}\right), 0 \right\}$ with the water temperature $W_s$ at time $s$ and a small parameter $\varepsilon = 10^{-4}$ to maintain the positivity of $w'_s$, and $\left(\bar{W}, \underline{W}\right) = (25, 5)$ (°C). This coefficient $4\left(\bar{W} - \underline{W}\right)^{-2}$ is a normalization factor. The water temperature $W_s$ as a sinusoidal function of time is identified using a least-squares fitting at Station Y (**Figure 5**):

$$W_s = 14.36 - 7.70 \cos\left( \frac{2\pi s}{P} \right) - 4.00 \sin\left( \frac{2\pi s}{P} \right) \ \ (\text{°C}), \ \ s \in [0, P], \tag{107}$$

where the time 0 corresponds to the beginning of January 1 of a year. In this baseline case, the modeled water temperature varies between 5 °C and 25 °C (maximum: 23.04 °C and minimum: 5.68 °C). However, as the water temperature may shift due to climate change in the future, we consider hypothetical scenarios with higher and lower temperatures as well.

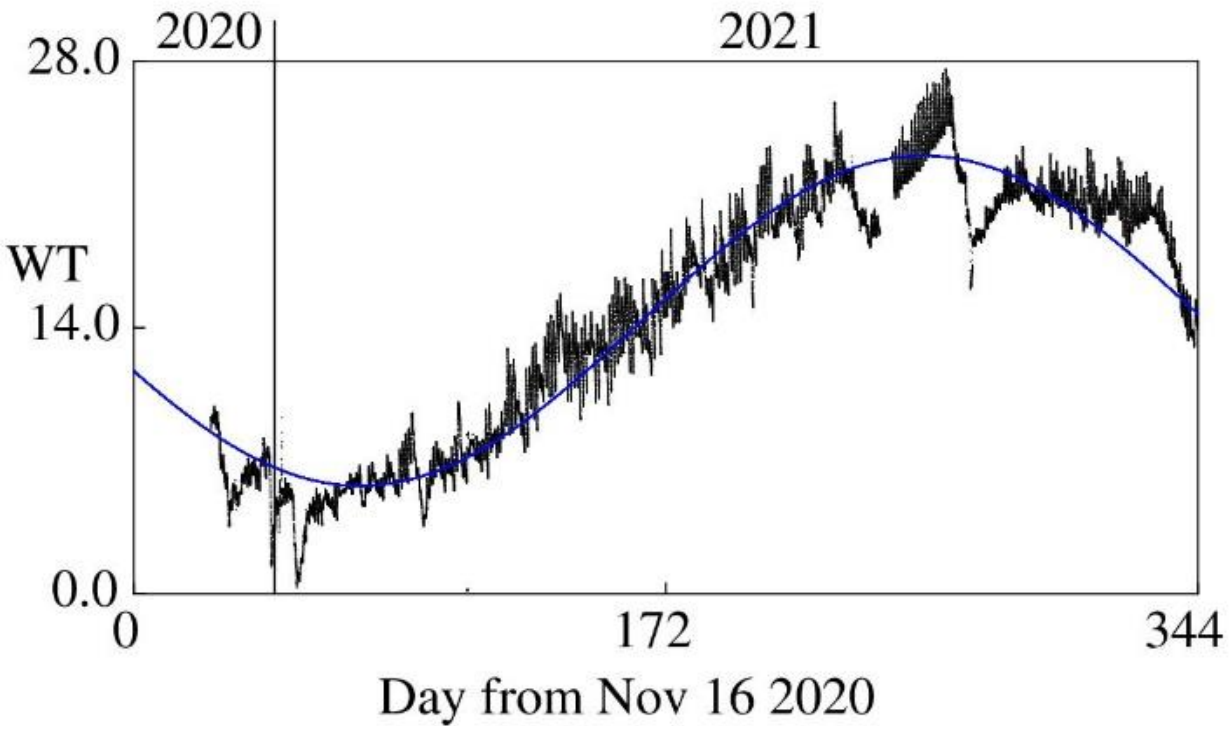

**Figure 5.** Observed (black) and fitted (blue) water temperature at Station Y.

#### 4.5.2 Demonstrative example

**Figures 6** and **7** show the demonstrative computational examples of the Riccati equation at Station Y with $w = 100$, suggesting that the quantities $\Gamma$ and $\gamma$ modulating the optimal control are seasonally varying. Especially, **Figure 7** shows that the optimal level of the discharge is controlled according to the weighting factor in a way that a relatively small discharge is preferred during winter (near the time 0 (year) and 1 (year) in the figure), while a relatively large discharge should be maintained in summer (near the time 0.6 (year) in the figure) during which the growth of benthic algae should be inspected carefully.

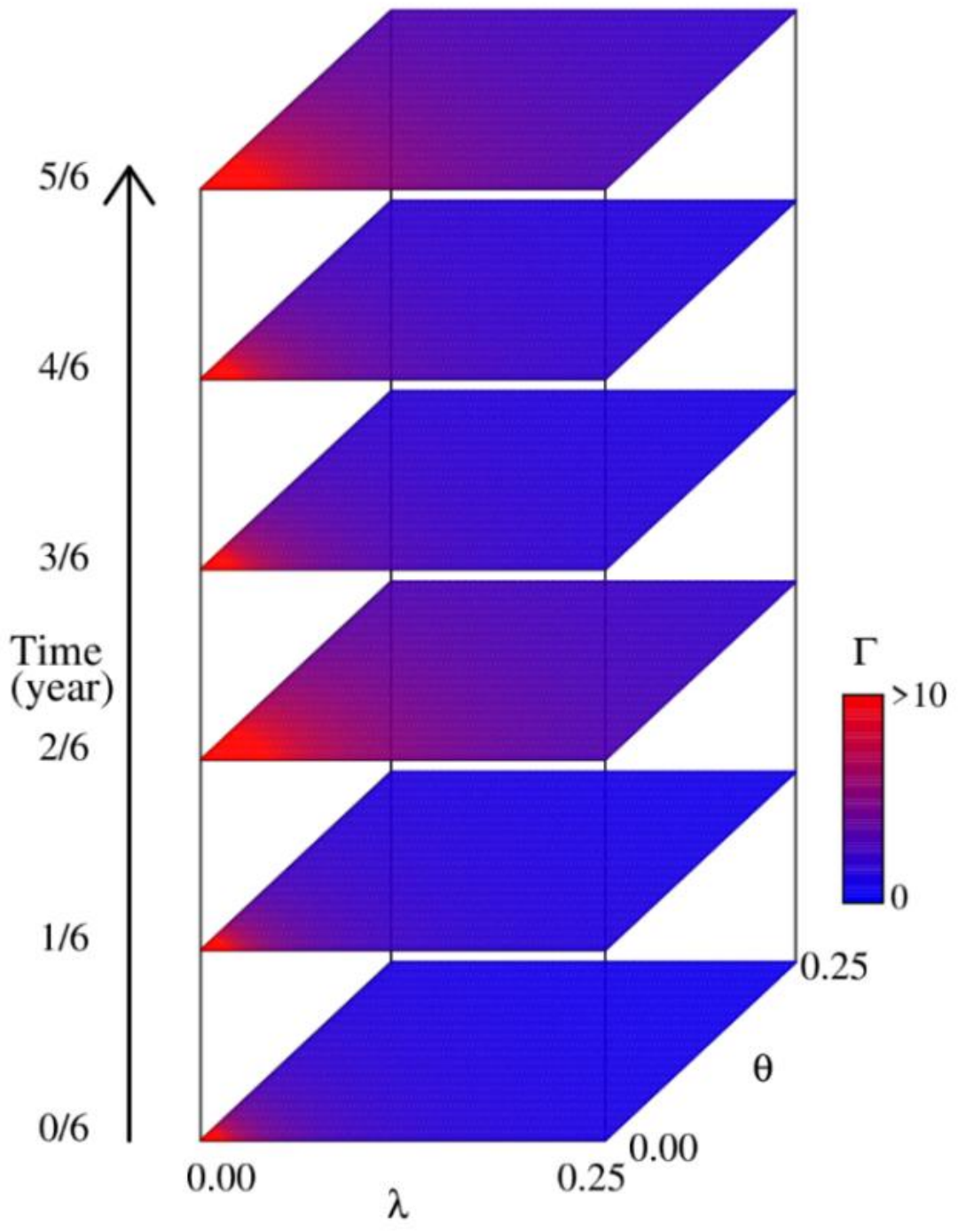

**Figure 6.** Computational example of $\Gamma$ at Station Y.

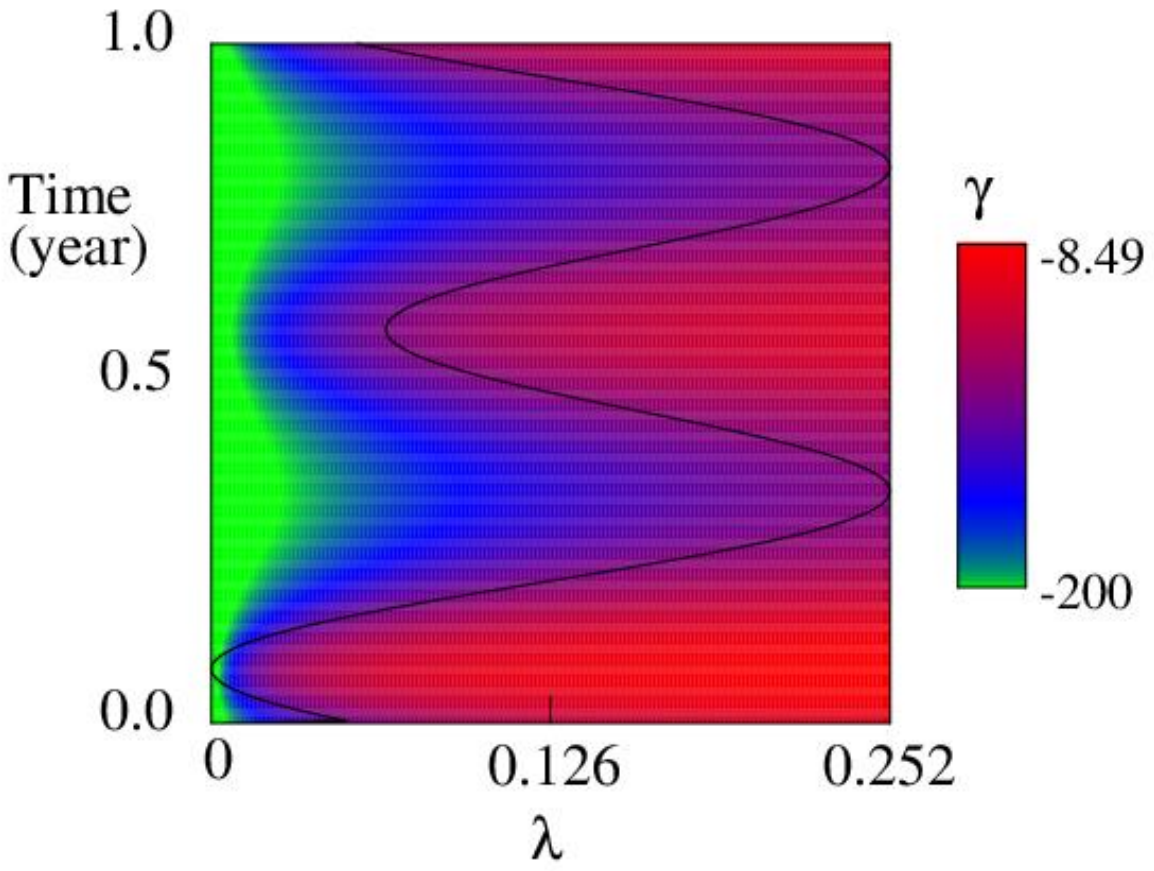

**Figure 7.** Computational example of $\gamma$ at Station Y. The weighting factor $w'$ is also plotted in the figure as a blackline normalized to fit the plot area.

#### 4.5.3 Cost and performance

We compute efficient frontiers between the cost and performance to analyze a continuum of each optimal control from the viewpoint of cost and benefit. Given a control $u$, set

$$D = \limsup_{k\to+\infty} \frac{1}{kP} \mathbb{E}\left[\int_0^{kP} \frac{w'_s}{2}\left(X_{n,s} - \hat{X}_s\right)^2 \mathrm{d}s\right] \text{ and } C = \limsup_{k\to+\infty} \frac{1}{kP} \mathbb{E}\left[\int_0^{kP} \frac{u_s^2}{2} \mathrm{d}s\right], \quad (108)$$

with which we recover $J = D + wC$ as a scalarized objective between the deviation $D$ and the control cost $C$ through a weighting factor $w$. The index $D$ is understood as a metric of stabilizability of the controlled dynamics because it measures probabilistic closeness between the controlled state and its target. We can analyze the relationship between cost and performance, namely stability, by concurrently considering $C$ and $D$.

Using optimal control $u = u^*$, we can compute an efficient frontier in the 2-D $C$-$D$ plane with $w > 0$. More specifically, given a sequence $\{w_i\}_{i=1,2,3,...}$, for each $w_i$, we solve the Riccati equation and obtain $J|_{w=w_i}$ to obtain the control $u^*|_{w=w_i}$, and then, we find $C|_{w=w_i}$ from the KBE. After that, we obtain $D|_{w=w_i} = J|_{w=w_i} - w_i C|_{w=w_i}$. Finally, we obtain points $(C,D)|_{w=w_i}$ on a frontier.

**Figure 8** shows the computed efficient frontier at each station. The deviation divided by the cost $D/C$ decreases in the order of Station Y, Station U, and Station D. This can be explained by the order of standard deviation decrease presented in **Table 2**. The deviation will be larger for an environment that has more significantly varying discharge time series. The stabilizability of the optimal control is supported in **Figure 8** because $D$ decreases as $C$ increases. This suggests that the supOU process can be stabilized so that the controlled state remains close to the target with high probability, even if the drift of the controlled dynamics is not necessarily dissipative. The computational results suggest the existence of a stabilization mechanism in the proposed control problem.

All the efficient frontiers are convex. Therefore, one can choose the best control in each frontier by imposing another criterion. For example, in practice, we can choose the one where the deviation $D$ is comparable to the variance of the uncontrolled state, namely the performance guarantee $D = c\mathrm{Std}^2$ with $c > 0$. This choice is reasonable in our case because $\mathrm{Std}^2$ has the order of $O(10^2)$ ($m^6/s^2$) in all the stations according to **Table 2**, while **Figure 8** suggests that it is possible to significantly reduce the deviation to be smaller than the uncontrolled $\mathrm{Std}^2$ values if a sufficient amount of effort is paid for controlling the discharge. The horizontal line of the performance guarantee $D = 0.05\mathrm{Std}^2$ is plotted in **Figure 8** for each station. In this case, the required controlling costs $C$ are 4.7 (Station Y), 1.1 (Station D), and 1.8 (Station U). One can choose the crossing point to single out the point of each efficient frontier.

Finally, we apply the model to hypothetical scenarios with lower and higher water temperatures at Station Y. We examine the proposed model with the water temperature (107) increased or decreased by 1, 2, and 3 (°C) to analyze the impact of the warming and cooling of the water temperature at Station Y on the cost and performance of optimal controls. The ranges of the increase/decrease of the river water temperature are reasonable in Japan based on recent investigations (Morid et al., 2020; Ye and Kameyama, 2021).

**Figures 9** and **10** show the computed efficient frontiers at Station Y for increased and decreased water temperatures, respectively. Again, all the computed efficient frontiers are convex. The horizontal line $D = 0.05\mathrm{Std}^2$ is also plotted in these figures as a demonstration. Both figures show that increasing and decreasing the water temperature at Station Y shifts the frontier downward, indicating that the deviation divided by the cost decreases. Therefore, climate changes leading to the temperature shift improve the efficiency of the optimal control of the discharge at Station Y. Indeed, the required controlling cost $C$ to achieve $D = 0.05\mathrm{Std}^2$ increases as the temperature regime deviates from the baseline.

The analysis here is simplified in the sense that it does not consider the dependence of watershed hydrological processes on climate change. Nevertheless, the proposed model can be applied to such cases if the temperature dependence of the supOU process of discharge is clarified. This is currently unresolved but will be addressed in the future in several rivers.

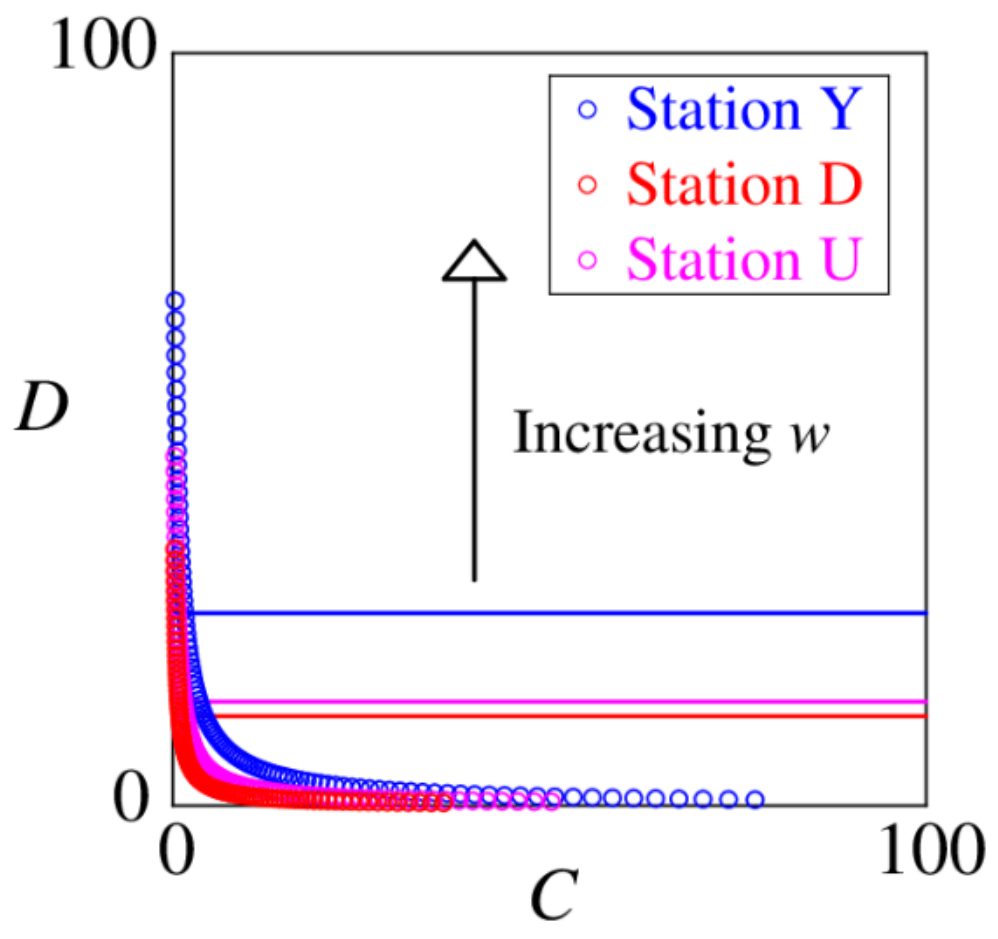

**Figure 8.** Computed efficient frontiers. The dimension of $C$ is the square of that of discharge: ($m^6/s^2$), while that of $D$ is the square of the acceleration of discharge: ($m^6/s^4$). Thus, the dimension of $w$ is ($s^2$). The parameter values $w = w_i$ are $w_i = 10^{-2+i/25}$ ($i = 0,1,2,...,100$). The horizontal lines $D = 0.05\mathrm{Std}^2$ are plotted with the corresponding colors. In this case, the required controlling costs $C$ are 2.3 (Station Y), 1.1 (Station D), and 1.8 (Station U).

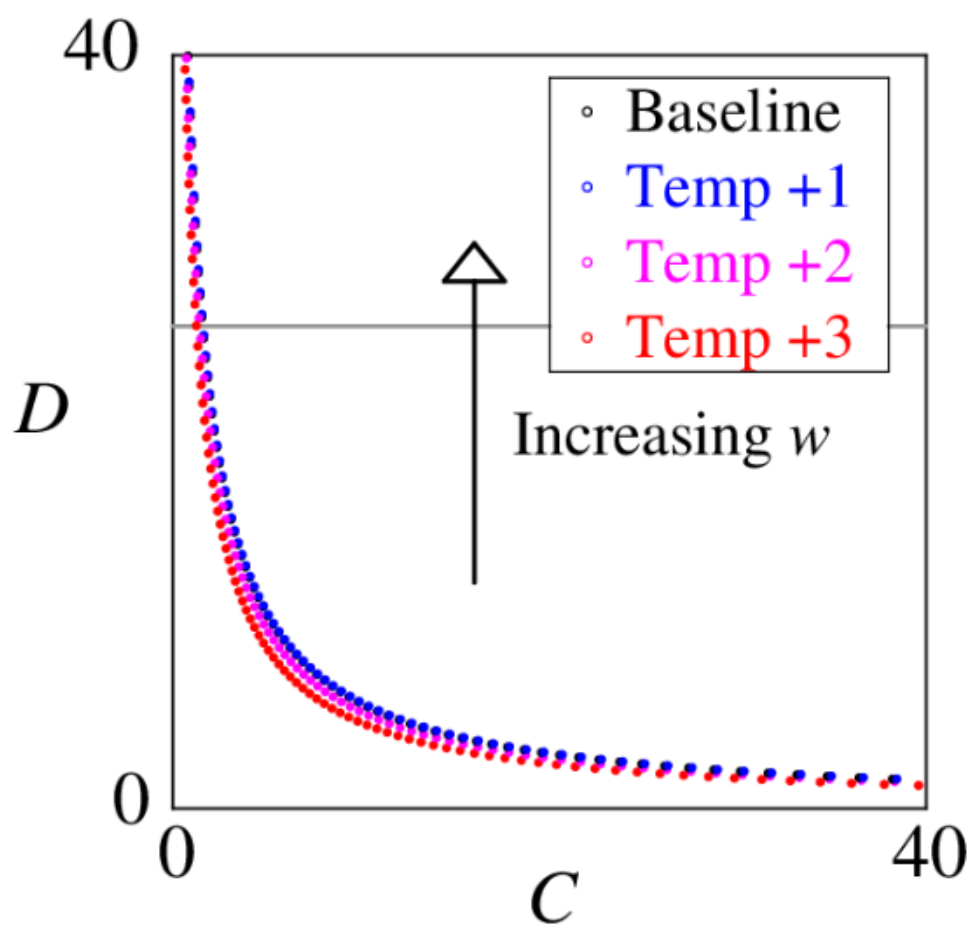

**Figure 9.** Computed efficient frontiers at Station Y for different water temperature scenarios. “Baseline” represents the case considered in (107), and “Temp +X” means the scenario with the water temperature increased by X (°C). The horizontal line $D = 0.05\mathrm{Std}^2$ is plotted in gray. In this case, the required controlling costs $C$ are 2.3 (Nominal), 2.3 (Temp +1), 2.1 (Temp +2), and 1.7 (Temp +3).

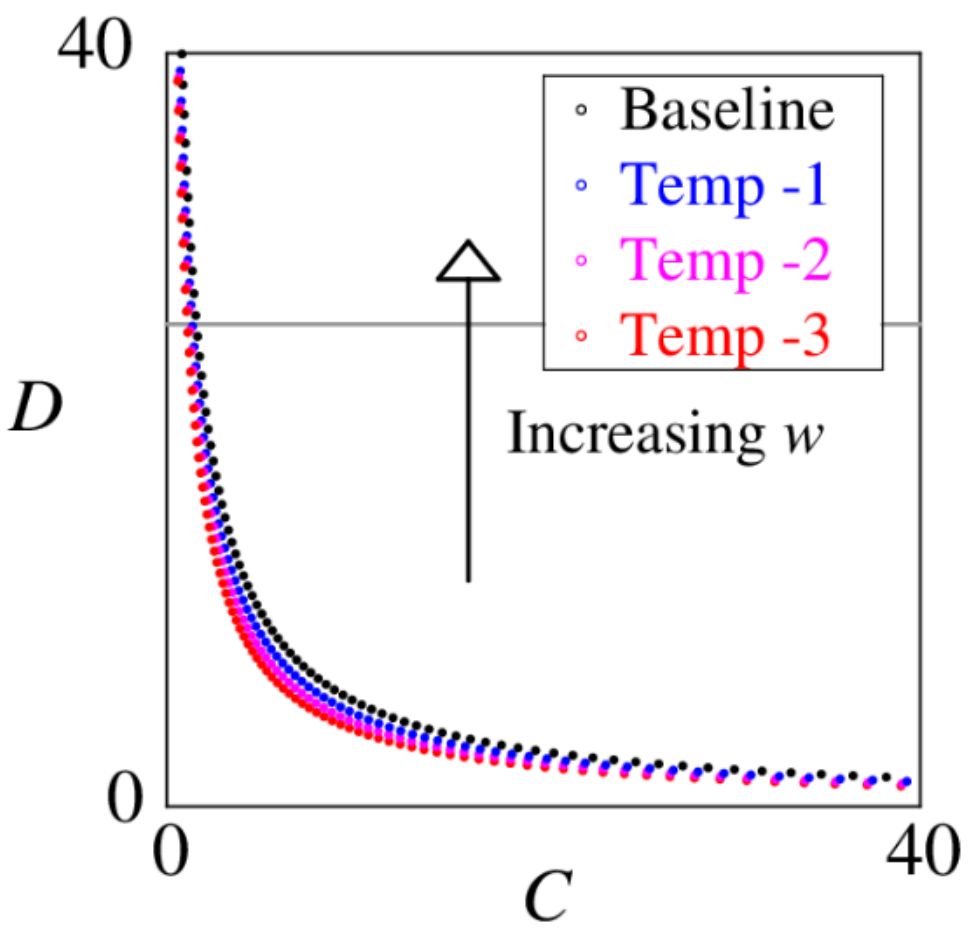

**Figure 10.** Computed efficient frontiers at Station Y for different water temperature scenarios. "Baseline" represents the case considered in (107), and "Temp −X" means the scenario with the water temperature decreased by X (°C). The horizontal line $D = 0.05\mathrm{Std}^2$ is plotted in gray. In this case, the required controlling costs $C$ are 1.6 (Temp −3), 1.8 (Temp −2), 2.0 (Temp −1), and 2.3 (Nominal).

## 5 Conclusion

We demonstrated that supOU processes driven by generalized tempered stable noises can be effectively used for the analysis of discharge time series data. The integral operator Riccati equation to determine the optimal control was derived from the infinite-dimensional HJB equation and its optimality was verified. In addition, the KBE was presented for evaluating the performance of optimal controls in a theoretically consistent and computable framework. The proposed model was applied to the real data. A manufactured solution was used to demonstrate the convergence of numerical solutions. An efficient frontier analysis revealed the balance between cost and performance.

The control of supOU processes is still at an evolving stage. The theoretical numerical analysis of the Riccati equation is unresolved, although our computational results are supportive. Mild solution approaches may open a door to a better understanding of the control problems of supOU processes (Peszat and Zabczyk, 2007). The applicability of higher-order numerical schemes, such as spectral schemes, will be a pivotal issue toward a more efficient approximation of the Riccati equation. Larger-scale LQ problems can be handled by sparse grids (Harbrecht and Kalmykov, 2021). Our approach for the performance evaluation via a KBE can be extended to other controls with slight modifications. Efficient solvers for Riccati and related equations are available, but their applicability to supOU cases is unknown (Hernández-Verón and Romero, 2018; Zhang et al., 2019), and will be examined in future.

Our problem is a specific LQ-type. More realistic problems would involve state constraints (Calvia et al., 2021) with which the HJB equation does not reduce to a Riccati equation. Solutions in such cases will be understood as viscosity solutions (Fabbri et al., 2017). The origin and impacts of non-smoothness should then be discussed because these characteristics can provide vital information about controlled dynamics. Coupling the discharge with sediment transport (Yoshioka et al., 2021) would lead to non-smooth dynamics. An efficient numerical scheme such as a tensor decomposition (Dolgov et al., 2021) will be necessary for dealing with non-LQ problems. Control delay may also exist in real cases. Another future direction will be considering model ambiguity based on an LQ framework (Brock et al., 2014). Currently, we are tackling this topic based on a robust control methodology combined with an LQ control formulation. Finally, considering a computable state-constraint approach would be interesting. Such a method exists for deterministic systems (Aubin-Frankowski, 2021), but its expendability to stochastic systems is an open question.

**Appendix A**

The hydraulic model (Tanaka et al., 2021) to generate the discharge time series at Station Y is briefly explained. It contains two sub-models: a 1-km mesh distributed hydrological sub-model (1K-DHM, Tanaka and Tachikawa, 2015) and an inundation model coupling rainfall–runoff models (IMCR, Tanaka et al., 2017), both of which have already been verified against real problems. 1K-DHM receives effective rainfall in the target catchment and simulates mountainous surface runoff and streamflow discharge based on a kinematic wave model. The simulated surface runoff is used as the input of IMCR to simulate downscaled river flows using 1-D local inertial equations with the spatial resolution ranging from 50 to 100 (m). The discharge time series at Station Y is the computed discharge of ICMR at a corresponding grid point.

The upstream and downstream boundary conditions of the IMCR were provided by the discharge data at Obara Dam (Station D) and the water level at Station H in **Figure 1**, respectively. Rainfall data were obtained from gauged rainfall data (blue circles in **Figure 1**) around the river and interpolated by the Thiessen polygon method. The effective rainfall of the catchment was estimated as rainfall minus evapotranspiration, which was estimated monthly using the common Thornthwaite method with the monthly mean temperature at the nearest meteorological gauge (green circle in **Figure 1**). The model parameters of 1K-DHM were calibrated with a severe flood event in July 2018, as explained in Tanaka et al. (2021) using the data at Station O. We validate this hydraulic model by comparing the simulated and observed water depths at Station Y from March to October, 2021 (**Figure A1**). They agree well at both low and high flows, except for the period after the second huge flood in August 2021. This is considered due to a flood-induced displacement of the riverbed on which the data logger was placed.

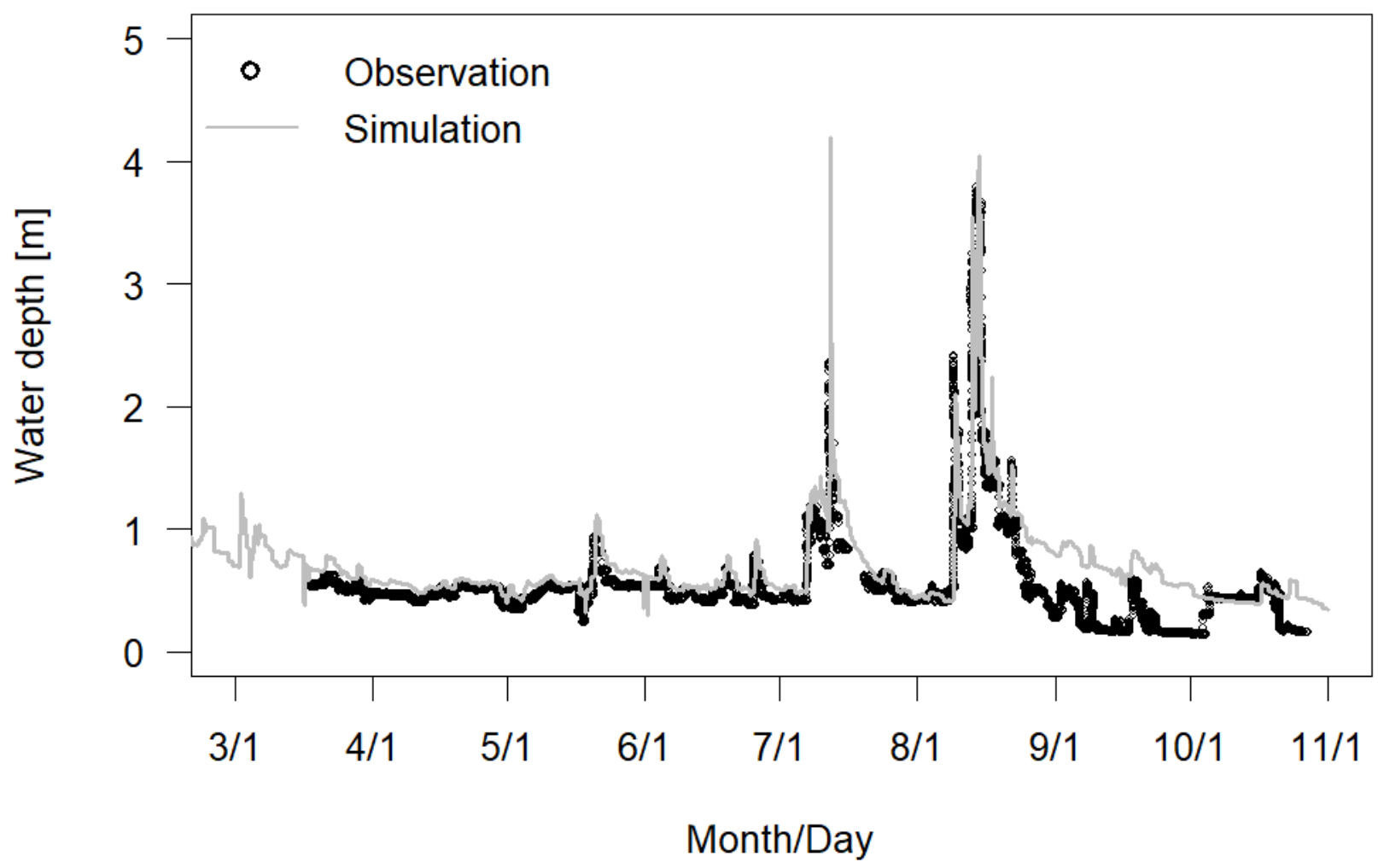

**Figure A1.** Observed (circles) and simulated (gray line) water depths at Station Y in 2021.

**Appendix B**

The identified tempered stable model as a simpler version of the generalized tempered stable model in the main text is presented in **Tables B1** and **B2**. The tempered stable model corresponds to the case $p_v = 1$ and hence has a weaker exponential decay of large jumps. **Table B2** shows that the minimized objective (95) at each observation station is larger than that presented in **Table 2**, supporting the use of the generalized tempered stable model.

**Table B1.** Identified parameter values.

| Parameter | Station Y | Station D | Station U |
|---|---|---|---|
| $B_\pi$ (1/h) | 0.0344 | 0.0201 | 0.0315 |
| $\alpha_\pi$ (-) | 2.17 | 2.97 | 2.53 |
| $\underline{X}$ (m$^3$/s) | 1.28 | 1.00 | 0.00 |
| $p_v$ (-) | 1 | 1 | 1 |
| $\alpha_v$ (-) | 6.68.E-01 | 1.02.E-01 | 4.601.E-01 |
| $a_v$ (m$^{3\alpha_v}$/s$^{\alpha_v}$) | 2.43.E-02 | 7.49.E-03 | 5.176.E-03 |
| $b_v$ (s/m$^3$) | 3.43.E-03 | 1.85.E-03 | 9.155.E-03 |

**Table B2.** Statistics of the identified model and empirical ones with their relative errors.

| | Station Y | | Station D | | Station U | |
|---|---|---|---|---|---|---|
| Least-squares error | 0.0141 | | 0.00885 | | 0.00517 | |
| | Data | Model | Data | Model | Data | Model |
| Ave (m$^3$/s) | 1.21.E+01 | 1.20.E+01 | 5.13.E+00 | 5.06.E+00 | 5.40.E+00 | 5.36.E+00 |
| Std (m$^3$/s) | 2.26.E+01 | 2.28.E+01 | 1.54.E+01 | 1.56.E+01 | 1.66.E+01 | 1.67.E+01 |
| Skew (-) | 1.28.E+01 | 1.14.E+01 | 1.19.E+01 | 1.08.E+01 | 1.27.E+01 | 1.19.E+01 |
| Kurt (-) | 2.43.E+02 | 2.54.E+02 | 1.95.E+02 | 2.01.E+02 | 2.54.E+02 | 2.61.E+02 |

## Appendix C

### C.1 Proof of Proposition 1

We compare the two characteristic functions. For $t,u \in \mathbb{R}$, we get

$$
\begin{aligned}
&\left|\ln\varphi_{X_t}(u)-\ln\varphi_{n,X_{n,t}}(u)\right| \\
&=\left|\begin{array}{l}\int_0^{+\infty}\int_{-\infty}^{t}\int_0^{+\infty}\left(\cos\left(ue^{-\lambda(t-s)}z\right)+\mathrm{i}\sin\left(ue^{-\lambda(t-s)}z\right)-1\right)v(\mathrm{d}z)\mathrm{d}s\pi(\mathrm{d}\lambda) \\ -\int_0^{+\infty}\int_{-\infty}^{t}\int_0^{+\infty}\left(\cos\left(ue^{-\lambda(t-s)}z\right)+\mathrm{i}\sin\left(ue^{-\lambda(t-s)}z\right)-1\right)v(\mathrm{d}z)\mathrm{d}s\pi_n(\mathrm{d}\lambda)\end{array}\right| \\
&\le\left|\begin{array}{l}\int_0^{+\infty}\int_{-\infty}^{t}\int_0^{+\infty}\left(\cos\left(ue^{-\lambda(t-s)}z\right)-1\right)v(\mathrm{d}z)\mathrm{d}s\pi(\mathrm{d}\lambda) \\ -\int_0^{+\infty}\int_{-\infty}^{t}\int_0^{+\infty}\left(\cos\left(ue^{-\lambda(t-s)}z\right)-1\right)v(\mathrm{d}z)\mathrm{d}s\pi_n(\mathrm{d}\lambda)\end{array}\right| \\
&\quad+\left|\begin{array}{l}\int_0^{+\infty}\int_{-\infty}^{t}\int_0^{+\infty}\sin\left(ue^{-\lambda(t-s)}z\right)v(\mathrm{d}z)\mathrm{d}s\pi(\mathrm{d}\lambda) \\ -\int_0^{+\infty}\int_{-\infty}^{t}\int_0^{+\infty}\sin\left(ue^{-\lambda(t-s)}z\right)v(\mathrm{d}z)\mathrm{d}s\pi_n(\mathrm{d}\lambda)\end{array}\right|
\end{aligned}, \tag{109}
$$

where $\exp\left(\mathrm{i}ue^{-\lambda(t-s)}z\right)=\cos\left(ue^{-\lambda(t-s)}z\right)+\mathrm{i}\sin\left(ue^{-\lambda(t-s)}z\right)$. The first term in the last line of (109) is handled as follows (we use the classical Fubini theorem to exchange the order of integrations by $\left|\cos\left(ue^{-\lambda(t-s)}z\right)-1\right|\le 2$ and $\left|\cos\left(ue^{-\lambda(t-s)}z\right)-1\right|\le|u|z$):

$$
\begin{aligned}
&\left|\begin{array}{l}\int_0^{+\infty}\int_{-\infty}^{t}\int_0^{+\infty}\left(\cos\left(ue^{-\lambda(t-s)}z\right)-1\right)v(\mathrm{d}z)\mathrm{d}s\pi(\mathrm{d}\lambda) \\ -\int_0^{+\infty}\int_{-\infty}^{t}\int_0^{+\infty}\left(\cos\left(ue^{-\lambda(t-s)}z\right)-1\right)v(\mathrm{d}z)\mathrm{d}s\pi_n(\mathrm{d}\lambda)\end{array}\right| \\
&=\left|\int_0^{+\infty}\int_{-\infty}^{t}\left\{\int_0^{+\infty}\left(\cos\left(ue^{-\lambda(t-s)}z\right)-1\right)\pi(\mathrm{d}\lambda)-\int_0^{+\infty}\left(\cos\left(ue^{-\lambda(t-s)}z\right)-1\right)\pi_n(\mathrm{d}\lambda)\right\}\mathrm{d}sv(\mathrm{d}z)\right| \\
&\le\int_0^{+\infty}\int_{-\infty}^{t}\left|\int_0^{+\infty}\left(\cos\left(ue^{-\lambda(t-s)}z\right)-1\right)\pi(\mathrm{d}\lambda)-\int_0^{+\infty}\left(\cos\left(ue^{-\lambda(t-s)}z\right)-1\right)\pi_n(\mathrm{d}\lambda)\right|\mathrm{d}sv(\mathrm{d}z)
\end{aligned}. \tag{110}
$$

In the rest of this sub-section, a generic positive constant independent from $i,n,t,u$ whose value may be different among different lines is denoted as $C_0$. We obtain

$$
\begin{aligned}
&\left|\int_0^{+\infty}\left(\cos\left(ue^{-\lambda(t-s)}z\right)-1\right)\pi(\mathrm{d}\lambda)-\int_0^{+\infty}\left(\cos\left(ue^{-\lambda(t-s)}z\right)-1\right)\pi_n(\mathrm{d}\lambda)\right| \\
&\le\left|\sum_{i=1}^{n}\left(\begin{array}{l}\int_{\eta_{n,i-1}}^{\eta_{n,i}}\left(\cos\left(ue^{-\lambda(t-s)}z\right)-1\right)\pi(\mathrm{d}\lambda) \\ -\int_{\eta_{n,i-1}}^{\eta_{n,i}}\left(\cos\left(ue^{-\lambda(t-s)}z\right)-1\right)\pi_n(\mathrm{d}\lambda)\end{array}\right)\right|+\int_{\eta_{n,n}}^{+\infty}\left(1-\cos\left(ue^{-\lambda(t-s)}z\right)\right)\pi(\mathrm{d}\lambda) \\
&\le\sum_{i=1}^{n}\left|\begin{array}{l}\int_{\eta_{n,i-1}}^{\eta_{n,i}}\left(\cos\left(ue^{-\lambda(t-s)}z\right)-1\right)\pi(\mathrm{d}\lambda) \\ -\int_{\eta_{n,i-1}}^{\eta_{n,i}}\left(\cos\left(ue^{-\lambda(t-s)}z\right)-1\right)\pi_n(\mathrm{d}\lambda)\end{array}\right|+\int_{\eta_{n,n}}^{+\infty}\left(1-\cos\left(ue^{-\lambda(t-s)}z\right)\right)\pi(\mathrm{d}\lambda)
\end{aligned}. \tag{111}
$$

We evaluate each term on the right-hand side. For its first term, we have

$$\begin{aligned}
&\left|\int_{\eta_{n,i-1}}^{\eta_{n,i}}\left(\cos\left(ue^{-\lambda(t-s)}z\right)-1\right)\pi(\mathrm{d}\lambda)-\int_{\eta_{n,i-1}}^{\eta_{n,i}}\left(\cos\left(ue^{-\lambda(t-s)}z\right)-1\right)\pi_n(\mathrm{d}\lambda)\right| \\
&=\left|\int_{\eta_{n,i-1}}^{\eta_{n,i}}\left(\cos\left(ue^{-\lambda(t-s)}z\right)-1\right)\pi(\mathrm{d}\lambda)-\left(\cos\left(ue^{-\lambda_i(t-s)}z\right)-1\right)c_i\right| \\
&=\left|\int_{\eta_{n,i-1}}^{\eta_{n,i}}\left(\cos\left(ue^{-\lambda(t-s)}z\right)-1\right)\pi(\mathrm{d}\lambda)-\int_{\eta_{n,i-1}}^{\eta_{n,i}}\left(\cos\left(ue^{-\lambda_i(t-s)}z\right)-1\right)\pi(\mathrm{d}\lambda)\right| \\
&=\left|\int_{\eta_{n,i-1}}^{\eta_{n,i}}\left(\cos\left(ue^{-\lambda(t-s)}z\right)-\cos\left(ue^{-\lambda_i(t-s)}z\right)\right)\pi(\mathrm{d}\lambda)\right| \\
&\leq|u|z\int_{\eta_{n,i-1}}^{\eta_{n,i}}\left|e^{-\lambda(t-s)}-e^{-\lambda_i(t-s)}\right|\pi(\mathrm{d}\lambda) \\
&=|u|z\left(\int_{\eta_{n,i-1}}^{\lambda_i}\left(e^{-\lambda(t-s)}-e^{-\lambda_i(t-s)}\right)\pi(\mathrm{d}\lambda)+\int_{\lambda_i}^{\eta_{n,i}}\left(e^{-\lambda_i(t-s)}-e^{-\lambda(t-s)}\right)\pi(\mathrm{d}\lambda)\right)
\end{aligned} \quad . \tag{112}$$

The last term of (111) is evaluated as

$$\int_{\eta_{n,n}}^{+\infty}\left(1-\cos\left(ue^{-\lambda(t-s)}z\right)\right)\pi(\mathrm{d}\lambda)\leq\int_{\eta_{n,n}}^{+\infty}|u|ze^{-\lambda(t-s)}\pi(\mathrm{d}\lambda)\leq|u|ze^{-\eta_{n,n}(t-s)}\int_{\eta_{n,n}}^{+\infty}\pi(\mathrm{d}\lambda)\,. \tag{113}$$

Using (112) and (113), we get

$$\begin{aligned}
&\left|\int_0^{+\infty}\left(\cos\left(ue^{-\lambda(t-s)}z\right)-1\right)\pi(\mathrm{d}\lambda)-\int_0^{+\infty}\left(\cos\left(ue^{-\lambda(t-s)}z\right)-1\right)\pi_n(\mathrm{d}\lambda)\right| \\
&\leq C_0|u|z\left\{\sum_{i=1}^{n}\left(\int_{\eta_{n,i-1}}^{\lambda_i}\left(e^{-\lambda(t-s)}-e^{-\lambda_i(t-s)}\right)\pi(\mathrm{d}\lambda)+\int_{\lambda_i}^{\eta_{n,i}}\left(e^{-\lambda_i(t-s)}-e^{-\lambda(t-s)}\right)\pi(\mathrm{d}\lambda)\right)+e^{-\eta_{n,n}(t-s)}\int_{\eta_{n,n}}^{+\infty}\pi(\mathrm{d}\lambda)\right\} \\
&\equiv C_0|u|zG(n,t-s)
\end{aligned} \; . \tag{114}$$

Similarly, we obtain

$$\left|\int_0^{+\infty}\sin\left(ue^{-\lambda(t-s)}z\right)\pi(\mathrm{d}\lambda)-\int_0^{+\infty}\sin\left(ue^{-\lambda(t-s)}z\right)\pi_n(\mathrm{d}\lambda)\right|\leq C_0|u|zG(n,t-s) \tag{115}$$

and thus

$$\left|\ln\varphi_{X_t}(u)-\ln\varphi_{n,X_{n,t}}(u)\right|\leq\left|\int_{-\infty}^{t}\int_0^{+\infty}C_0|u|zG(n,t-s)v(\mathrm{d}z)\mathrm{d}s\right|=C_0|u|M_1\int_{-\infty}^{t}G(n,t-s)\mathrm{d}s\,. \tag{116}$$

Now, we have

$$\begin{aligned}
&\int_{-\infty}^{t}\sum_{i=1}^{n}\left(\int_{\eta_{n,i-1}}^{\lambda_i}\left(e^{-\lambda(t-s)}-e^{-\lambda_i(t-s)}\right)\pi(\mathrm{d}\lambda)+\int_{\lambda_i}^{\eta_{n,i}}\left(e^{-\lambda_i(t-s)}-e^{-\lambda(t-s)}\right)\pi(\mathrm{d}\lambda)\right)\mathrm{d}s \\
&=\sum_{i=1}^{n}\int_{-\infty}^{t}\left(\int_{\eta_{n,i-1}}^{\lambda_i}\left(e^{-\lambda(t-s)}-e^{-\lambda_i(t-s)}\right)\pi(\mathrm{d}\lambda)+\int_{\lambda_i}^{\eta_{n,i}}\left(e^{-\lambda_i(t-s)}-e^{-\lambda(t-s)}\right)\pi(\mathrm{d}\lambda)\right)\mathrm{d}s \\
&=\sum_{i=1}^{n}\left(\int_{\eta_{n,i-1}}^{\lambda_i}\int_{-\infty}^{t}\left(e^{-\lambda(t-s)}-e^{-\lambda_i(t-s)}\right)\mathrm{d}s\,\pi(\mathrm{d}\lambda)+\int_{\lambda_i}^{\eta_{n,i}}\int_{-\infty}^{t}\left(e^{-\lambda_i(t-s)}-e^{-\lambda(t-s)}\right)\mathrm{d}s\,\pi(\mathrm{d}\lambda)\right)
\end{aligned} \tag{117}$$

because the absolute value of each integrand is bounded by 1. By the assumption of $\alpha_\pi > 2$, we can proceed as

$$
\begin{aligned}
&\sum_{i=1}^{n}\left(\int_{\eta_{n,i-1}}^{\lambda_i}\int_{-\infty}^{t}\left(e^{-\lambda(t-s)}-e^{-\lambda_i(t-s)}\right)\mathrm{d}s\pi(\mathrm{d}\lambda)+\int_{\lambda_i}^{\eta_{n,i}}\int_{-\infty}^{t}\left(e^{-\lambda_i(t-s)}-e^{-\lambda(t-s)}\right)\mathrm{d}s\pi(\mathrm{d}\lambda)\right)\\
&=\sum_{i=1}^{n}\left(\int_{\eta_{n,i-1}}^{\lambda_i}\left(\frac{1}{\lambda}-\frac{1}{\lambda_i}\right)\pi(\mathrm{d}\lambda)+\int_{\lambda_i}^{\eta_{n,i}}\left(\frac{1}{\lambda_i}-\frac{1}{\lambda}\right)\pi(\mathrm{d}\lambda)\right)\\
&=\sum_{i=1}^{n}\frac{1}{\lambda_i}\left(\int_{\eta_{n,i-1}}^{\lambda_i}\frac{\lambda_i-\lambda}{\lambda}\pi(\mathrm{d}\lambda)+\int_{\lambda_i}^{\eta_{n,i}}\frac{\lambda-\lambda_i}{\lambda}\pi(\mathrm{d}\lambda)\right)\\
&\le\sum_{i=1}^{n}\frac{1}{\lambda_i}\int_{\eta_{n,i-1}}^{\eta_{n,i}}\frac{|\lambda_i-\lambda|}{\lambda}\pi(\mathrm{d}\lambda)\\
&\le\sum_{i=1}^{n}\frac{\eta_{n,i}-\eta_{n,i-1}}{\lambda_i}\int_{\eta_{n,i-1}}^{\eta_{n,i}}\frac{1}{\lambda}\pi(\mathrm{d}\lambda)\\
&\le C_0\sum_{i=1}^{n}\frac{\left(\eta_{n,i}-\eta_{n,i-1}\right)^2}{\lambda_i}
\end{aligned}
\quad , \qquad (118)
$$

where we used

$$
\int_{\eta_{n,i-1}}^{\eta_{n,i}}\frac{1}{\lambda}\pi(\mathrm{d}\lambda)\le C_0\int_{\eta_{n,i-1}}^{\eta_{n,i}}\lambda^{\alpha_\pi-2}e^{-\frac{\lambda}{B_\pi}}\mathrm{d}\lambda\le C_0\left(\eta_{n,i}-\eta_{n,i}\right) \qquad (119)
$$

as $\alpha_\pi>2$ and thus the integrand is uniformly bounded. Hence, we obtain

$$
\begin{aligned}
\left|\ln\varphi_{X_t}(u)-\ln\varphi_{n,X_{n,t}}(u)\right|&\le C_0|u|M_1\int_{-\infty}^{t}G(n,t-s)\mathrm{d}s\\
&\le C_0|u|M_1\left\{\sum_{i=1}^{n}\frac{\left(\eta_{n,i}-\eta_{n,i-1}\right)^2}{\lambda_i}+\frac{1}{\eta_{n,n}}\int_{\eta_{n,n}}^{+\infty}\pi(\mathrm{d}\lambda)\right\}.\\
&=C_0|u|M_1\bar{G}(n)
\end{aligned}
\qquad (120)
$$

Because (120) holds for any $t,u\in\mathbb{R}$, we obtain $\left|\ln\varphi_{X_t}(u)-\ln\varphi_{n,X_{n,t}}(u)\right|\to 0$ as $n\to+\infty$ if $\bar{G}(n)\to+0$ under the same limit. Note that $\left|\varphi_{X_t}(u)-\varphi_{n,X_{n,t}}(u)\right|\to 0$ as $n\to+\infty$ if $\left|\ln\varphi_{X_t}(u)-\ln\varphi_{n,X_{n,t}}(u)\right|\to 0$ under this limit because $\varphi_{X_t}(u)>0$ is finite. The proof is completed by (116)–(120) along with the uniqueness and convergence results of the characteristic functions (Theorem 1.4, Theorem 3.1 of Vandegriffe (2020)).

### C.2 Proof of Proposition 2

We rewrite (42) as (recall that is $A$ symmetric)

$$
\begin{aligned}
0&=\frac{\partial A_{i,j}}{\partial s}-\left(\lambda_i+\lambda_j\right)A_{i,j}-\frac{1}{w}\sum_{k=1}^{n}c_kA_{k,i}\sum_{k=1}^{n}c_kA_{k,j}+1\\
&=\frac{\partial A_{i,j}}{\partial s}-\left(\lambda_i+\lambda_j\right)A_{i,j}-\sum_{k,l=1}^{n}A_{i,k}\frac{1}{\sqrt{w}}c_k\frac{1}{\sqrt{w}}c_lA_{l,j}+1 \qquad \left(1\le i,j\le n,\ 0\le s\le P\right) \qquad (121)\\
&=\frac{\partial A_{i,j}}{\partial s}+\sum_{k=1}^{n}\bar{A}_{i,k}A_{k,j}+\sum_{k=1}^{n}A_{i,k}\bar{A}_{j,k}-\sum_{k,l,m=1}^{n}A_{i,k}\bar{C}_{m,k}\bar{C}_{m,l}A_{l,j}+\sum_{k=1}^{n}\bar{B}_{i,k}\bar{B}_{j,k}
\end{aligned}
$$

with

$$\bar{A} = \left[\bar{A}_{i,j}\right]_{1\le i,j\le n} = \left[-\delta_{i,j}\lambda_i\right]_{1\le i,j\le n}, \quad \bar{B} = \left[\bar{B}_{i,j}\right]_{1\le i,j\le n} = \left[\delta_{i,j}\right]_{1\le i,j\le n}, \quad \bar{C} = \left[\bar{C}_{i,j}\right]_{1\le i,j\le n} = \left[\frac{1}{\sqrt{w}}c_j\right]_{1\le i,j\le n}. \tag{122}$$

The conditions to be checked to apply to Theorem 6 of Ha and Choi (2019) are the stabilizability of the pair $\left(\bar{A},\bar{B}\right)$ and the detectability of the pair $\left(\bar{A},\bar{C}\right)$.

The stabilizability is satisfied if $\left(\bar{A},\bar{B}\right)$ is controllable. It holds if the rank of the following matrix $m_{1,n}$ equals $n$:

$$m_{1,n} = \left[\bar{B},\bar{A}\bar{B},\bar{A}^2\bar{B},...,\bar{A}^{n-1}\bar{B}\right], \text{ where } \bar{A}^k\bar{B} = \left[\sum_{m=1}^{n}\left(\delta_{i,m}\left(-\lambda_i\right)^k\delta_{m,j}\right)\right]_{1\le i,j\le n} = \left[\left(-\lambda_i\right)^k\right]_{1\le i,j\le n}. \tag{123}$$

By $0<\lambda_1<\lambda_2,...,<\lambda_n$, a straightforward matrix transformation gives

$$\operatorname{rank}\left(\bar{A}^k\bar{B}\right) = \operatorname{rank}\left[\left(\lambda_i\right)^k\right]_{1\le i,j\le n} \text{ and hence } \operatorname{rank}\left(m_{1,n}\right) = \operatorname{rank}\left[\left(\lambda_i\right)^{j-1}\right]_{1\le i,j\le n}. \tag{124}$$

This $\left[\left(\lambda_i\right)^{j-1}\right]_{1\le i,j\le n}$ is a Vandermonde matrix (Liesen and Mehrmann, 2015), showing $\operatorname{rank}\left(m_{1,n}\right)=n$.

The detectability follows if the rank of the following matrix $m_{2,n}$ equals $n$:

$$m_{2,n} = \left[\bar{C},\bar{A}\bar{C},\bar{A}^2\bar{C},...,\bar{A}^{n-1}\bar{C}\right], \quad \bar{A}^k\bar{C} = \left[\sum_{m=1}^{n}\left(\delta_{i,m}\left(-\lambda_i\right)^k\frac{1}{\sqrt{w}}c_j\right)\right]_{1\le i,j\le n} = \left[\frac{1}{\sqrt{w}}\left(-\lambda_i\right)^k c_j\right]_{1\le i,j\le n}. \tag{125}$$

By $0<\lambda_1<\lambda_2,...,<\lambda_n$ and $c_i>0$ ($i=1,2,...,n$), a matrix transformation yields

$$\operatorname{rank}\left(\bar{A}^k\bar{C}\right) = \bar{A}^k\bar{C} = \operatorname{rank}\left[\left(-\lambda_i\right)^k c_j\right]_{1\le i,j\le n} = \operatorname{rank}\left[\left(-\lambda_i\right)^k\right]_{1\le i,j\le n} = \operatorname{rank}\left[\left(\lambda_i\right)^k\right]_{1\le i,j\le n}. \tag{126}$$

We again obtain the Vandermonde matrix: $\operatorname{rank}\left(m_{2,n}\right)=n$. Consequently, (42) admits a unique positive semi-definite $A$. Given this $A$, we find a solution to the linear system (43). Because its right-hand side is strictly decreasing and Lipschitz continuous with respect to $B_i$ ($i=1,2,...,n$), its unique solvability in a point-wise sense follows Proposition 1 and Theorem 1 and Bostan (2002). The uniqueness of $H_n$ is trivial.

### C.3 Proof of Proposition 3

Set a $u\in\mathbb{U}$ and a corresponding controlled process $Y$ with $Y_t\in L^2\left(\pi\right)$ a.s. $t>0$. By, Itô's formula (Theorem D.2 of Peszat and Zabczyk, 2007), we have

$$
\begin{aligned}
&\mathrm{d}\left(\frac{1}{2}\int_0^{+\infty}\int_0^{+\infty}\Gamma_s(\theta,\lambda)Y_s(\lambda)Y_s(\theta)\pi(\mathrm{d}\lambda)\pi(\mathrm{d}\theta)+\int_0^{+\infty}\gamma_s(\lambda)Y_s(\lambda)\pi(\mathrm{d}\lambda)\right)\\
&=\frac{1}{2}\int_0^{+\infty}\int_0^{+\infty}\mathrm{d}\Gamma_s(\theta,\lambda)Y_s(\lambda)Y_s(\theta)\pi(\mathrm{d}\lambda)\pi(\mathrm{d}\theta)+\int_0^{+\infty}\mathrm{d}\gamma_s(\lambda)Y_s(\lambda)\pi(\mathrm{d}\lambda)\\
&+\frac{1}{2}\int_0^{+\infty}\int_0^{+\infty}\Gamma_s(\theta,\lambda)\mathrm{d}\left(Y_s(\lambda)Y_s(\theta)\right)\pi(\mathrm{d}\lambda)\pi(\mathrm{d}\theta)+\int_0^{+\infty}\gamma_s(\lambda)\mathrm{d}Y_s(\lambda)\pi(\mathrm{d}\lambda)\\
&\equiv I_1+I_2+I_3+I_4
\end{aligned} \quad . \tag{127}
$$

By (67) and (68), we obtain

$$
\begin{aligned}
I_1&=\frac{1}{2}\int_0^{+\infty}\int_0^{+\infty}\mathrm{d}\Gamma_s(\theta,\lambda)Y_s(\lambda)Y_s(\theta)\pi(\mathrm{d}\lambda)\pi(\mathrm{d}\theta)\\
&=\frac{1}{2}\int_0^{+\infty}\int_0^{+\infty}\left(\begin{aligned}&(\theta+\lambda)\Gamma_s(\theta,\lambda)\\&+\frac{1}{w}\int_0^{+\infty}\Gamma_s(\theta_1,\lambda)\pi(\mathrm{d}\theta_1)\int_0^{+\infty}\Gamma_s(\theta_2,\theta)\pi(\mathrm{d}\theta_2)-1\end{aligned}\right)Y_s(\lambda)Y_s(\theta)\pi(\mathrm{d}\lambda)\pi(\mathrm{d}\theta)\mathrm{d}s
\end{aligned} \tag{128}
$$

and

$$
\begin{aligned}
I_2&=\int_0^{+\infty}\mathrm{d}\gamma_s(\lambda)Y_s(\lambda)\pi(\mathrm{d}\lambda)\\
&=\int_0^{+\infty}\left(\begin{aligned}&\lambda\gamma_s(\lambda)+\frac{1}{w}\int_0^{+\infty}\Gamma_s(\theta,\lambda)\pi(\mathrm{d}\lambda)\int_0^{+\infty}\gamma_s(\lambda)\pi(\mathrm{d}\lambda)\\&-M_1\int_0^{+\infty}\Gamma_s(\theta,\lambda)\pi(\mathrm{d}\lambda)+\bar{X}_s\end{aligned}\right)Y_s(\lambda)\pi(\mathrm{d}\lambda)\mathrm{d}s
\end{aligned} \quad . \tag{129}
$$

Furthermore, we have

$$
\begin{aligned}
I_4&=\int_0^{+\infty}\gamma_s(\lambda)\mathrm{d}Y_s(\lambda)\pi(\mathrm{d}\lambda)\\
&=\int_0^{+\infty}\left\{\left(-\lambda\gamma_s(\lambda)Y_s(\lambda)+u_s\gamma_s(\lambda)\right)\mathrm{d}s+\gamma_s\mathrm{d}L_s(\lambda)\right\}\pi(\mathrm{d}\lambda)
\end{aligned} \quad . \tag{130}
$$

The term $\mathrm{d}\left(Y_s(\lambda)Y_s(\theta)\right)$ is handled as follows (jumps have finite variations): if $\lambda\neq\theta$, then,

$$
\begin{aligned}
\mathrm{d}\left(Y_s(\lambda)Y_s(\theta)\right)&=\left(-\lambda Y_s(\lambda)+u_s\right)Y_s(\theta)\mathrm{d}s+Y_s(\theta)\mathrm{d}L_s(\lambda)\\
&+\left(-\theta Y_s(\theta)+u_s\right)Y_s(\lambda)\mathrm{d}s+Y_s(\lambda)\mathrm{d}L_s(\theta)\\
&=\left\{-(\lambda+\theta)Y_s(\lambda)Y_s(\theta)+u_s\left(Y_s(\theta)+Y_s(\lambda)\right)\right\}\mathrm{d}s\\
&+Y_s(\lambda)\mathrm{d}L_s(\theta)+Y_s(\theta)\mathrm{d}L_s(\lambda)
\end{aligned} \tag{131}
$$

while if $\lambda=\theta$, then,

$$
\begin{aligned}
d\left(Y_s(\lambda)Y_s(\theta)\right)&=\mathrm{d}\left(Y_s(\lambda)\right)^2\\
&=2\left\{-\lambda\left(Y_s(\lambda)\right)^2+u_sY_s(\lambda)\right\}\mathrm{d}s+\left(Y_{s-}(\lambda)+\mathrm{d}L_s(\lambda)\right)^2-\left(Y_{s-}(\lambda)\right)^2\\
&=2\left\{-\lambda\left(Y_s(\lambda)\right)^2+u_sY_s(\lambda)\right\}\mathrm{d}s+2Y_{s-}(\lambda)\mathrm{d}L_s(\lambda)+\left(\mathrm{d}L_s(\lambda)\right)^2
\end{aligned} \quad . \tag{132}
$$

Then, $I_3$ is rewritten as

$$\frac{1}{2}\int_0^{+\infty}\int_0^{+\infty}\Gamma_s(\theta,\lambda)\mathrm{d}\left(Y_s(\lambda)Y_s(\theta)\right)\pi(\mathrm{d}\lambda)\pi(\mathrm{d}\theta)$$

$$=\frac{1}{2}\int_0^{+\infty}\int_0^{+\infty}\Gamma_s(\theta,\lambda)\left\{\begin{array}{l}\left\{-(\lambda+\theta)Y_s(\lambda)Y_s(\theta)+u_s\left(Y_s(\theta)+Y_s(\lambda)\right)\right\}\mathrm{d}s\\+Y_s(\lambda)\mathrm{d}L_s(\theta)+Y_s(\theta)\mathrm{d}L_s(\lambda)\\+\chi_{\{\lambda=\theta\}}\left(\mathrm{d}L_s(\lambda)\right)^2\end{array}\right\}\pi(\mathrm{d}\lambda)\pi(\mathrm{d}\theta). \tag{133}$$

Consequently, we have

$$\begin{aligned}&I_1+I_2+I_3+I_4\\&=\frac{1}{2}\int_0^{+\infty}\int_0^{+\infty}\left(\begin{array}{l}(\theta+\lambda)\Gamma_s(\theta,\lambda)\\+\frac{1}{w}\int_0^{+\infty}\Gamma_s(\theta_1,\lambda)\pi(\mathrm{d}\theta_1)\int_0^{+\infty}\Gamma_s(\theta_2,\theta)\pi(\mathrm{d}\theta_2)-1\end{array}\right)Y_s(\lambda)Y_s(\theta)\pi(\mathrm{d}\lambda)\pi(\mathrm{d}\theta)\mathrm{d}s\\&+\int_0^{+\infty}\left(\begin{array}{l}\lambda\gamma_s(\lambda)+\frac{1}{w}\int_0^{+\infty}\Gamma_s(\theta,\lambda)\pi(\mathrm{d}\theta)\int_0^{+\infty}\gamma_s(\theta)\pi(\mathrm{d}\theta)\\-M_1\int_0^{+\infty}\Gamma_s(\theta,\lambda)\pi(\mathrm{d}\lambda)+\bar{X}\end{array}\right)Y_s(\lambda)\pi(\mathrm{d}\lambda)\mathrm{d}s\\&+\frac{1}{2}\int_0^{+\infty}\int_0^{+\infty}\Gamma_s(\theta,\lambda)\left\{\begin{array}{l}\left\{-(\lambda+\theta)Y_s(\lambda)Y_s(\theta)+u_s\left(Y_s(\theta)+Y_s(\lambda)\right)\right\}\mathrm{d}s\\+Y_s(\lambda)\mathrm{d}L_s(\theta)+Y_s(\theta)\mathrm{d}L_s(\lambda)\\+\chi_{\{\lambda=\theta\}}\left(\mathrm{d}L_s(\lambda)\right)^2\end{array}\right\}\pi(\mathrm{d}\lambda)\pi(\mathrm{d}\theta)\\&+\int_0^{+\infty}\gamma_s(\lambda)\left\{\left(-\lambda Y_s(\lambda)+u_s\right)\mathrm{d}s+\mathrm{d}L_s(\lambda)\right\}\pi(\mathrm{d}\lambda)\end{aligned}. \tag{134}$$

This is rearranged as

$$\begin{aligned}&\frac{1}{2}\int_0^{+\infty}\int_0^{+\infty}\left(\frac{1}{w}\int_0^{+\infty}\Gamma_s(\theta_1,\lambda)\pi(\mathrm{d}\theta_1)\int_0^{+\infty}\Gamma_s(\theta_2,\theta)\pi(\mathrm{d}\theta_2)-1\right)Y_s(\lambda)Y_s(\theta)\pi(\mathrm{d}\lambda)\pi(\mathrm{d}\theta)\mathrm{d}s\\&+\int_0^{+\infty}\left(\begin{array}{l}\frac{1}{w}\int_0^{+\infty}\Gamma_s(\theta,\lambda)\pi(\mathrm{d}\theta)\int_0^{+\infty}\gamma_s(\theta)\pi(\mathrm{d}\theta)\\-M_1\int_0^{+\infty}\Gamma_s(\theta,\lambda)\pi(\mathrm{d}\lambda)+\bar{X}_s\end{array}\right)Y_s(\lambda)\pi(\mathrm{d}\lambda)\mathrm{d}s\\&+\frac{1}{2}\int_0^{+\infty}\int_0^{+\infty}\Gamma_s(\theta,\lambda)\left\{\begin{array}{l}u_s\left(Y_s(\theta)+Y_s(\lambda)\right)\mathrm{d}s+2Y_s(\lambda)\mathrm{d}L_s(\theta)\\+\chi_{\{\lambda=\theta\}}\left(\mathrm{d}L_s(\lambda)\right)^2\end{array}\right\}\pi(\mathrm{d}\lambda)\pi(\mathrm{d}\theta)\\&+\int_0^{+\infty}\gamma_s(\lambda)\left\{u_s\mathrm{d}s+\mathrm{d}L_s(\lambda)\right\}\pi(\mathrm{d}\lambda)\\&=\frac{1}{w}\left\{\begin{array}{l}\frac{1}{2}\int_0^{+\infty}\int_0^{+\infty}\int_0^{+\infty}\Gamma_s(\theta_1,\lambda)\pi(\mathrm{d}\theta_1)\int_0^{+\infty}\Gamma_s(\theta_2,\theta)\pi(\mathrm{d}\theta_2)Y_s(\lambda)Y_s(\theta)\pi(\mathrm{d}\lambda)\pi(\mathrm{d}\theta)\mathrm{d}s\\+\int_0^{+\infty}\int_0^{+\infty}\Gamma_s(\theta,\lambda)\pi(\mathrm{d}\theta)\int_0^{+\infty}\gamma_s(\theta)\pi(\mathrm{d}\theta)Y_s(\lambda)\pi(\mathrm{d}\lambda)\mathrm{d}s\end{array}\right\}\\&-\frac{1}{2}\int_0^{+\infty}\int_0^{+\infty}Y_s(\lambda)Y_s(\theta)\pi(\mathrm{d}\lambda)\pi(\mathrm{d}\theta)\mathrm{d}s+\int_0^{+\infty}\bar{X}_sY_s(\lambda)\pi(\mathrm{d}\lambda)\mathrm{d}s\\&-M_1\int_0^{+\infty}\int_0^{+\infty}\Gamma_s(\theta,\lambda)\pi(\mathrm{d}\lambda)Y_s(\lambda)\pi(\mathrm{d}\lambda)\mathrm{d}s+\int_0^{+\infty}\int_0^{+\infty}\Gamma_s(\theta,\lambda)Y_s(\lambda)\mathrm{d}L_s(\theta)\pi(\mathrm{d}\lambda)\pi(\mathrm{d}\theta)\\&+u_s\left\{\int_0^{+\infty}\int_0^{+\infty}\Gamma_s(\theta,\lambda)Y_s(\theta)\pi(\mathrm{d}\lambda)\pi(\mathrm{d}\theta)\mathrm{d}s+\int_0^{+\infty}\gamma_s(\lambda)\pi(\mathrm{d}\lambda)\mathrm{d}s\right\}\\&+\int_0^{+\infty}\gamma_s(\lambda)\mathrm{d}L_s(\lambda)\pi(\mathrm{d}\lambda)+\frac{1}{2}\int_0^{+\infty}\int_0^{+\infty}\Gamma_s(\theta,\lambda)\chi_{\{\lambda=\theta\}}\left(\mathrm{d}L_s(\lambda)\right)^2\pi(\mathrm{d}\lambda)\pi(\mathrm{d}\theta)\end{aligned}. \tag{135}$$

Taking the expectation of (135) yields

$$
\begin{aligned}
&\frac{1}{w}\mathbb{E}\left[\begin{array}{l}\frac{1}{2}\int_0^{+\infty}\int_0^{+\infty}\int_0^{+\infty}\Gamma_s(\theta_1,\lambda)\pi(\mathrm{d}\theta_1)\int_0^{+\infty}\Gamma_s(\theta_2,\theta)\pi(\mathrm{d}\theta_2)Y_s(\lambda)Y_s(\theta)\pi(\mathrm{d}\lambda)\pi(\mathrm{d}\theta)\\ +\int_0^{+\infty}\int_0^{+\infty}\Gamma_s(\theta,\lambda)\pi(\mathrm{d}\theta)\int_0^{+\infty}\gamma_s(\theta)\pi(\mathrm{d}\theta)Y_s(\lambda)\pi(\mathrm{d}\lambda)\end{array}\right]\mathrm{d}s\\
&+\mathbb{E}\left[-\frac{1}{2}\int_0^{+\infty}\int_0^{+\infty}Y_s(\lambda)Y_s(\theta)\pi(\mathrm{d}\lambda)\pi(\mathrm{d}\theta)+\int_0^{+\infty}\bar{X}_sY_s(\lambda)\pi(\mathrm{d}\lambda)\right]\mathrm{d}s\\
&+\mathbb{E}\left[u_s\left\{\int_0^{+\infty}\int_0^{+\infty}\Gamma_s(\theta,\lambda)Y_s(\theta)\pi(\mathrm{d}\lambda)\pi(\mathrm{d}\theta)+\int_0^{+\infty}\gamma_s(\lambda)\pi(\mathrm{d}\lambda)\right\}\right]\mathrm{d}s\\
&+M_1\int_0^{+\infty}\gamma_s(\lambda)\pi(\mathrm{d}\lambda)\mathrm{d}s+\frac{1}{2}\int_0^{+\infty}\int_0^{+\infty}\Gamma_s(\theta,\lambda)\delta_{\{\lambda=\theta\}}\frac{\mathrm{d}\lambda}{\pi(\mathrm{d}\lambda)}M_2\pi(\mathrm{d}\lambda)\pi(\mathrm{d}\theta)\mathrm{d}s\\
&=\frac{1}{w}\mathbb{E}\left[\frac{1}{2}\left(\int_0^{+\infty}\nabla V(Y_s)\pi(\mathrm{d}\lambda)\right)^2-\frac{1}{2}\left(\int_0^{+\infty}\gamma_s(\lambda)\pi(\mathrm{d}\lambda)\right)^2+wu_s\int_0^{+\infty}\nabla V(Y_s)\pi(\mathrm{d}\lambda)\right]\mathrm{d}s\\
&+\mathbb{E}\left[\frac{(\bar{X}_s)^2}{2}-\frac{1}{2}\left(\int_0^{+\infty}Y_s(\lambda)\pi(\mathrm{d}\lambda)-\bar{X}\right)^2\right]\mathrm{d}s\\
&+M_1\int_0^{+\infty}\gamma_s(\lambda)\pi(\mathrm{d}\lambda)\mathrm{d}s+\frac{1}{2}M_2\int_0^{+\infty}\Gamma_s(\lambda,\lambda)\pi(\mathrm{d}\lambda)\mathrm{d}s\\
&=\mathbb{E}\left[\frac{1}{2w}\left(\int_0^{+\infty}\nabla V(Y_s)\pi(\mathrm{d}\lambda)+wu_s\right)^2-\frac{1}{2}\left(\int_0^{+\infty}Y_s(\lambda)\pi(\mathrm{d}\lambda)-\bar{X}\right)^2-\frac{w}{2}u_s^2\right]\mathrm{d}s\\
&+M_1\int_0^{+\infty}\gamma_s(\lambda)\pi(\mathrm{d}\lambda)\mathrm{d}s+\frac{1}{2}M_2\int_0^{+\infty}\Gamma_s(\lambda,\lambda)\pi(\mathrm{d}\lambda)\mathrm{d}s\\
&+\frac{(\bar{X}_s)^2}{2}\mathrm{d}s-\frac{1}{2w}\left(\int_0^{+\infty}\gamma_s(\lambda)\pi(\mathrm{d}\lambda)\right)^2\mathrm{d}s
\end{aligned}
\quad . \tag{136}
$$

We integrate (136) with respect to $s$ from 0 to $kP$ ($k\in\mathbb{N}$), and use (68) and (69) to obtain

$$
\begin{aligned}
&\left[\mathbb{E}\left[\frac{1}{2}\int_0^{+\infty}\int_0^{+\infty}\Gamma_s(\theta,\lambda)Y_s(\lambda)Y_s(\theta)\pi(\mathrm{d}\lambda)\pi(\mathrm{d}\theta)+\int_0^{+\infty}\gamma_s(\lambda)Y_s(\lambda)\pi(\mathrm{d}\lambda)\right]\right]_{s=0}^{s=kP}\\
&=\frac{1}{2w}\mathbb{E}\left[\int_0^{kP}\left(\int_0^{+\infty}\nabla V(Y_s)\pi(\mathrm{d}\lambda)+wu_s\right)^2\mathrm{d}s\right]\\
&-\mathbb{E}\left[\int_0^{kP}\left(\frac{1}{2}\left(\int_0^{+\infty}Y_s(\lambda)\pi(\mathrm{d}\lambda)-\bar{X}_s\right)^2+\frac{w}{2}u_s^2\right)\mathrm{d}s\right]\\
&+\int_0^{kP}\left\{\begin{array}{l}M_1\int_0^{+\infty}\gamma_s(\lambda)\pi(\mathrm{d}\lambda)+\frac{1}{2}M_2\int_0^{+\infty}\Gamma_s(\lambda,\lambda)\pi(\mathrm{d}\lambda)\\ -\frac{1}{2w}\left(\int_0^{+\infty}\gamma_s(\lambda)\pi(\mathrm{d}\lambda)\right)^2+\frac{(\bar{X}_s)^2}{2}\end{array}\right\}\mathrm{d}s\\
&=kPH-\mathbb{E}\left[\int_0^{kP}\left(\frac{1}{2}\left(\int_0^{+\infty}Y_s(\lambda)\pi(\mathrm{d}\lambda)-\bar{X}_s\right)^2+\frac{w}{2}u_s^2\right)\mathrm{d}s\right]\\
&+\frac{1}{2w}\mathbb{E}\left[\int_0^{kP}\left(\int_0^{+\infty}\nabla V(Y_s)\pi(\mathrm{d}\lambda)+wu_s\right)^2\mathrm{d}s\right]
\end{aligned}
\quad . \tag{137}
$$

Through periodicity, we obtain

$$
\begin{aligned}
&\frac{1}{kP}\mathbb{E}\left[\int_0^{kP}\left(\frac{1}{2}\left(\int_0^{+\infty}Y_s(\lambda)\pi(\mathrm{d}\lambda)-\bar{X}_s\right)^2+\frac{w}{2}u_s^2\right)\mathrm{d}s\right] \\
&=H+\frac{1}{2wkP}\mathbb{E}\left[\int_0^{kP}\left(\int_0^{+\infty}\nabla V(Y_s)\pi(\mathrm{d}\lambda)+wu_s\right)^2\mathrm{d}s\right]
\end{aligned}. \tag{138}
$$

We then get

$$
J(u)=H+\frac{1}{2w}\limsup_{k\to+\infty}\frac{1}{kP}\mathbb{E}\left[\int_0^{kP}\left(\frac{1}{w}\int_0^{+\infty}\nabla V(Y_s)\pi(\mathrm{d}\lambda)+u_s\right)^2\mathrm{d}s\right]. \tag{139}
$$

The right-hand side is minimized by (56), and the minimum equals $H$. Consequently, we obtain the proposition.

**Appendix D**

We consider a deterministic and time-independent case, which allows us to analytically study the controlled dynamics. Now, the integral operator Riccati equation reads

$$0=-(\theta+\lambda)\Gamma(\theta,\lambda)-\frac{1}{w}\int_0^{+\infty}\Gamma(\theta_1,\lambda)\pi(\mathrm{d}\theta_1)\int_0^{+\infty}\Gamma(\theta_2,\theta)\pi(\mathrm{d}\theta_2)+1, \tag{140}$$

$$0=-\lambda\gamma(\lambda)-\frac{1}{w}\int_0^{+\infty}\Gamma(\theta,\lambda)\pi(\mathrm{d}\theta)\int_0^{+\infty}\gamma(\lambda)\pi(\mathrm{d}\lambda)-\bar{X}, \tag{141}$$

$$H=-\frac{1}{2w}\left(\int_0^{+\infty}\gamma(\lambda)\pi(\mathrm{d}\lambda)\right)^2+\frac{1}{2}\bar{X}^2. \tag{142}$$

We rearrange (140) as

$$\begin{aligned}&\left(\frac{1}{\theta}+\frac{1}{\lambda}\right)\Gamma(\theta,\lambda)=-\frac{1}{w}\left(\frac{1}{\lambda}\int_0^{+\infty}\Gamma(\theta_1,\lambda)\pi(\mathrm{d}\theta_1)\right)\left(\frac{1}{\theta}\int_0^{+\infty}\Gamma(\theta_2,\theta)\pi(\mathrm{d}\theta_2)\right)+\frac{1}{\theta}\frac{1}{\lambda}\\&\int_0^{+\infty}\int_0^{+\infty}\left(\frac{1}{\theta}+\frac{1}{\lambda}\right)\Gamma(\theta,\lambda)\pi(\mathrm{d}\lambda)\pi(\mathrm{d}\theta)\\&=-\frac{1}{w}\left(\int_0^{+\infty}\int_0^{+\infty}\frac{1}{\lambda}\Gamma(\theta_1,\lambda)\pi(\mathrm{d}\theta_1)\pi(\mathrm{d}\lambda)\right)\left(\int_0^{+\infty}\int_0^{+\infty}\frac{1}{\theta}\Gamma(\theta_2,\theta)\pi(\mathrm{d}\theta_2)\pi(\mathrm{d}\theta)\right)\\&+\left(\int_0^{+\infty}\frac{1}{\theta}\pi(\mathrm{d}\theta)\right)\left(\int_0^{+\infty}\frac{1}{\lambda}\pi(\mathrm{d}\lambda)\right)\end{aligned}, \tag{143}$$

so

$$\frac{I^2}{w}=-2I+R^2 \quad \text{with} \quad I=\int_0^{+\infty}\int_0^{+\infty}\frac{1}{\lambda}\Gamma(\theta,\lambda)\pi(\mathrm{d}\lambda)\pi(\mathrm{d}\theta) \tag{144}$$

or more explicitly,

$$\frac{I}{\sqrt{w}}=-\sqrt{w}+\sqrt{w+R} \tag{145}$$

by the positive semi-definiteness of $\Gamma$. Similarly, we rearrange (141) as

$$\int_0^{+\infty}\gamma(\omega)\pi(\mathrm{d}\omega)=-w\bar{X}R\left(w+\int_0^{+\infty}\int_0^{+\infty}\frac{1}{\lambda}\Gamma(\theta,\lambda)\pi(\mathrm{d}\theta)\pi(\mathrm{d}\lambda)\right)^{-1} \tag{146}$$

and obtain

$$\int_0^{+\infty}\gamma(\omega)\pi(\mathrm{d}\omega)=-\frac{w}{w+I}R\bar{X}. \tag{147}$$

Under the control $u=u^*$, for $\lambda>0$ we have

$$\mathrm{d}Y_t(\lambda)=\left(-\lambda Y_t(\lambda)-\frac{1}{w}\left(\int_0^{+\infty}\int_0^{+\infty}\Gamma(\theta,\omega)Y_s(\theta)\pi(\mathrm{d}\theta)\pi(\mathrm{d}\omega)+\int_0^{+\infty}\gamma(\theta)\pi(\mathrm{d}\theta)\right)\right)\mathrm{d}t. \tag{148}$$

Assuming a stationary state obtains

$$Y_\infty(\lambda)=-\frac{1}{w\lambda}\left(\int_0^{+\infty}\int_0^{+\infty}\Gamma(\theta,\omega)Y_\infty(\theta)\pi(\mathrm{d}\theta)\pi(\mathrm{d}\omega)+\int_0^{+\infty}\gamma(\theta)\pi(\mathrm{d}\theta)\right), \tag{149}$$

leading to

$$\int_0^{\infty}\int_0^{+\infty}\Gamma(\theta,\lambda)Y_\infty(\lambda)\pi(\mathrm{d}\lambda)\pi(\mathrm{d}\theta)=-\frac{I}{w}\left(\int_0^{+\infty}\int_0^{+\infty}\Gamma(\theta,\omega)Y_\infty(\theta)\pi(\mathrm{d}\theta)\pi(\mathrm{d}\omega)+\int_0^{+\infty}\gamma(\theta)\pi(\mathrm{d}\theta)\right) \tag{150}$$

and hence

$$\int_0^{\infty}\int_0^{+\infty}\Gamma(\theta,\lambda)Y_{\infty}(\lambda)\pi(\mathrm{d}\lambda)\pi(\mathrm{d}\theta)=-\frac{I}{w+I}\int_0^{+\infty}\gamma(\theta)\pi(\mathrm{d}\theta). \tag{151}$$

By (149), we obtain

$$Y_{\infty}(\lambda)=-\frac{1}{w\lambda}\left(-\frac{I}{w+I}\int_0^{+\infty}\gamma(\theta)\pi(\mathrm{d}\theta)+\int_0^{+\infty}\gamma(\theta)\pi(\mathrm{d}\theta)\right)=-\frac{1}{\lambda(w+I)}\int_0^{+\infty}\gamma(\theta)\pi(\mathrm{d}\theta) \tag{152}$$

and the corresponding stationary state

$$X_{\infty}=\underline{X}-\frac{R}{w+I}\int_0^{+\infty}\gamma(\theta)\pi(\mathrm{d}\theta)=\underline{X}+\frac{1}{\left(\sqrt{w}+I/\sqrt{w}\right)^2}R^2\overline{X}\ . \tag{153}$$

From (144), under the no-cost limit $w\to+0$, (145) gives

$$\lim_{w\to+0}\frac{I}{\sqrt{w}}=\lim_{w\to+0}\left(-\sqrt{w}+\sqrt{w+R^2}\right)=R \tag{154}$$

and

$$X_{\infty}=\underline{X}+\frac{1}{(0+R)^2}R^2\left(\hat{X}-\underline{X}\right)=\hat{X}\ , \tag{155}$$

Then, the target $\hat{X}$ is an equilibrium of the controlled system. Moreover, (142) gives

$$H=-\frac{w}{2(w+I)^2}R^2\overline{X}^2+\frac{1}{2}\overline{X}^2\to 0 \ \text{ as } \ w\to 0\,, \tag{156}$$

suggesting that the state eventually approaches to $X=\hat{X}$ under $w\to+0$. We observe a similar behavior in the periodically driven stochastic case (**Section 4.5.3**).

Finally, the stationary state may be different from $\hat{X}$ due to the finite-dimensional assumption of the control (also see **Remark 3**). Nevertheless, the stationary state can be made arbitrary closer to the target by choosing a sufficiently small $w$.

**Declarations of interest:** There are no conflicts of interest to declare.

**Funding:** This work was supported by the Japan Society for the Promotion of Science [grant number 19H03073]; Kurita Water and Environment Foundation [grant number 21K008]; and Environmental Research Projects from the Sumitomo Foundation [grant number 203160]. Part of the computational work was conducted using the supercomputer of ACCMS, Kyoto University, Japan.